\documentclass[11pt]{article}
\usepackage[margin=1in]{geometry}
\usepackage{amsmath}
\usepackage{authblk}
\usepackage{amssymb}
\usepackage{amsthm}
\usepackage{color}
\usepackage[english]{babel}
\usepackage{mathtools}
\usepackage{cite}
\usepackage{caption}
\usepackage{hyperref}
\usepackage{multirow}
\usepackage{graphicx} 
\usepackage{algorithm}
\usepackage{algpseudocode}
\usepackage{epstopdf}
\usepackage{url}
\newtheorem{theorem}{Theorem}

\newtheorem{proposition}[theorem]{Proposition}
\newtheorem{corollary}[theorem]{Corollary}
\theoremstyle{remark}
\newtheorem{remark}[theorem]{Remark}

\newcommand{\R}{\mathbb{R}}

\renewcommand{\P}{\mathbb{P}}

\author{Christopher Kennedy\thanks{Email: \texttt{ckennedy@math.utexas.edu}.} }

\author{Rachel Ward\thanks{Email: \texttt{rward@math.utexas.edu}.}}

\affil{Department of Mathematics, University of Texas at Austin}

\title{Greedy Variance Estimation for the LASSO}

\begin{document}

 \maketitle
\begin{abstract}
Recent results have proven the minimax optimality of LASSO and related algorithms for noisy linear regression.  However, these results tend to rely on variance estimators that are inefficient or optimizations that are slower than LASSO itself.  We propose an efficient estimator for the noise variance in high dimensional linear regression that is faster than LASSO, only requiring $p$ matrix-vector multiplications.  We prove this estimator is consistent with a good rate of convergence, under the condition that the design matrix satisfies the Restricted Isometry Property (RIP).  In practice, our estimator scales incredibly well into high dimensions, is highly parallelizable, and only incurs a modest bias.
\end{abstract}

\section{Introduction}
The LASSO \cite{tibshirani1996regression} is a classical algorithm for doing noisy linear regression in the case when the number of regression coefficients is larger than the number of response variables $p > n$. The analysis of LASSO has recently surged with much work on establishing oracle inequalities for $\ell_2$ estimation over sparse vectors, and corresponding minimax rates.  Typically, such results rely on knowledge of the variance of the noise, which is unknown in practice.  The full extent of the literature on LASSO is immense and beyond the scope of this paper, but we point to a few important references on the oracle inequalities and corresponding minimax error rates (see \cite{bickel2009simultaneous}, \cite{meinshausen2009lasso}, \cite{zhang2008sparsity}, \cite{van2008high}, \cite{raskutti2011minimax}, \cite{zhang2009some}, \cite{wainwright2009information}, \cite{lounici2011oracle}, \cite{ye2010rate}, \cite{verzelen2012minimax}, \cite{candes2013well}).

A good review of variance estimators for LASSO is given in \cite{reid2013study}, where variance estimation using cross-validated LASSO is highlighted as particularly strong in many sparsity regimes.  This method typically uses 5 or 10-fold cross-validation to train the hyperparameters in LASSO and analysis relies on the restricted eigenvalue condition on the design matrix.  The above work was later complemented by a theoretical analysis of a slightly modified variant of cross-validated LASSO in \cite{chatterjee2015prediction} (see also \cite{fan2012variance} \cite{homrighausen2013lasso}, e.g.).  The method of moments (see \cite{dicker2014variance}) is a reasonable alternative to cross-validated LASSO.  It relies on the assumption that the design matrix is Gaussian and exploits statistical properties to formulate an estimator.  It is consistent with a good rate of convergence \cite{dicker2014variance}, but the design matrix has to be Gaussian which is restrictive.  We should also mention a variant of the LASSO - the square-root LASSO (see \cite{belloni2011square}) - whose penalty level doesn't depend on the variance of the noise.  However, the resulting estimator is formulated as a conic programming problem which can be inefficient in practice and is beyond the scope of this work.

\subsection{Our Contribution}
The main contributions of our a paper are the following:
\begin{itemize}
    \item We provide a \emph{fast} variance estimator.  In fact, our variance estimator only requires $p$ matrix-vector multiplications.  This also ensures our method is highly parallelizable, and is faster than a single iteration of LASSO.
    \item Our estimator is consistent in the sense that it converges in probability as $n, p \rightarrow \infty$ to the true variance, under mild conditions on the asymptotic growth of $\| \beta \|_2$.  We have a quantitative bound on the rate of convergence, as discussed following Theorem \ref{Theorem_orth_mat}.
    \item We only require a deterministic assumption on the design matrix (the Restricted Isometry Property) which holds with high probability over many standard matrix ensembles over appropriate parameter regimes.  In particular, the condition holds for any orthonormal design matrix.
\end{itemize}

Our estimator admits a surprisingly simple theoretical argument for convergence using standard compressed sensing-type results and concentration of measure.  We note that we assume the design matrix satisfies the restricted isometry property, which is stronger than the restricted eigenvalue condition typically considered in the literature (see \cite{van2009conditions} e.g.) but still admits a wide range of random matrix ensembles.

In practice, although our estimator exhibits a higher bias than more standard estimators like cross-validated LASSO, it does well in the high-dimensional regime where the method of moments estimator become prohibitive to compute.

\subsection{Notation}

For a matrix $X \in \R^{n\times p}$, and a subset $\Omega\subset \{1,..,p\}$, $X_{\Omega} \in \R^{n \times |\Omega|}$ will denote the restriction of $X$ to its columns indexed by $\Omega$.\\
For a vector $v \in \R^p$, $\Omega_v$ is defined to be the support of $v$.\\
For each $j = 1,...,p/L$, we use $\Omega_j := \{(j-1)L+1,...,jL\}$ to denote the $j$th ``window" of the signal.

\section{Problem Statement}

Suppose $\beta \in \mathbb{R}^p$ is $s$-sparse, and that we are given a noisy, transformed version of this signal:
$$
y= X\beta + \eta,
$$
where $\eta \in \mathbb{R}^n$ has i.i.d. Gaussian entries $\eta_j \sim {\cal N}(0, \sigma^2)$ and $X \in \R^{n\times p}$ is a known design matrix.  For the purpose of analysis, we define a notion of a well-behaved design matrix $X$.  We will assume that the matrix $X$ satisfies the Restricted Isometry Property (RIP), which was introduced in \cite{candes2005decoding} and is a common property used in Compressed Sensing.  It guarantees that a matrix is a near-isometry on sparse vectors.  Specifically, we say $X$ satisfies the RIP of order $s$ and level $\delta>0$ if for all $z$ such that $\|z\|_0 \leq s$,
$$
(1-\delta)\|z\|_2^2 \leq \| Xz \|_2^2 \leq (1 + \delta)\|z\|_2^2.
$$
RIP of order $s$ and level $\delta > 0$ is satisfied with probability at least $1- 2\exp(-\delta^2 n/(2C))$ once
\begin{equation}
n \geq C \delta^{-2} s \log(p/s)
\label{s}
\end{equation}
on an $n \times p$ matrix $X$ whose entries $X_{i,j}$ are independent realizations of a subgaussian random variable with mean 0 and variance $1/n$, such as a Gaussian or Bernoulli random variable \cite{foucart2013mathematical, baraniuk2008simple}.  Here, $C > 0$ is a universal constant independent of all parameters. The Restricted Isometry Property is obtained with high probability on many classes of structured random matrices, such as random partial Fourier matrices \cite{rudelson2008sparse, rauhut2010compressive}, but with a slightly smaller (by factors of $\log(p)$) constant $s$.

It is in general nontrivial to recover the true signal $\beta$.  We consider here the standard LASSO algorithm \cite{tibshirani1996regression} to return a denoised version of $\beta$:

\begin{align}
\label{LASSO_denoise}
\widehat{\beta} &= \text{arg} \min_\beta   \| X\beta - y \|_2^2 + 2 \lambda \| \beta \|_1 .
\end{align}
The magnitude of the parameter $\lambda$ in the objective (\ref{LASSO_denoise}) controls the balance between the $\ell_1$ term which promotes sparsity in the recovered signal $\widehat{\beta},$ and a mean squared error term $\| X \beta - y\|_2^2$ which promotes consistency with the observed measurements.  It is important to balance the two terms appropriately so that one doesn't overfit to the transformed signal $y$ but also doesn't over-enforce sparsity.  The standard analysis of the LASSO is conditioned on the event $\{ \lambda: \lambda/4 \geq \| X^T \eta\|_{\infty} /n\}$ (see \cite{bickel2009simultaneous}).  In particular, for the case that $\eta$ is Gaussian with variance $\sigma^2$ and $X$ is orthogonal, with high probability we have $\|X^T\|_{\infty}/n = \Theta(\sigma^2 \log(n) /n)$.  Thus, with the choice $\lambda = 4 \sigma^2 \log(n)/n$, the LASSO will provably produce a good estimate $\beta$.

However, in applications, the variance $\sigma$, and hence a proper choice of $\lambda$, is not known a priori.  We consider the case where $\sigma$ is not known in advance, and needs to be estimated from the signal $y$.  It should be clear from the above observations that precision in estimating the parameter $\sigma$ improves recovery of the true signal.

\section{Greedy Variance Estimation -- The Orthonormal Case}
\label{sect:ortho_design}
For the moment we focus on the case where $X\in \R^{p\times p}$ is an orthonormal matrix ($p=n$) and the problem reduces to recovering the noisy signal $y = \beta + \eta$ (by rotational invariance of the Gaussian).
In this regime, the LASSO has the closed form solution
$$
\widehat{\beta}_i = \text{sign}(y_i)(|y_i| - \lambda)_+,
$$
where $\widehat{\beta_i} = \widehat{\beta_i}(\lambda)$ implicitly depends on $\lambda$.  A standard approach is to minimize the cross-validation error: $$ \min_{\lambda} \|y - \widehat{\beta}(\lambda)\|_2,
$$
which has nice practical and theoretical properties (see \cite{kohavi1995study} e.g.).  Moreover, given the optimal $\lambda$ one can infer a good estimate of the variance as $\|\widehat{\beta} - y\|_2/p$.  However, this approach still requires one to compute the LASSO minimizer over a range of $\lambda$ values, whereas one would like to perform a single computation to estimate the variance (and thus optimal $\lambda$).  We formulate a method to estimate the variance which only needs a single pass over the input $y$.

\begin{figure}
\centering
\includegraphics[scale=0.25]{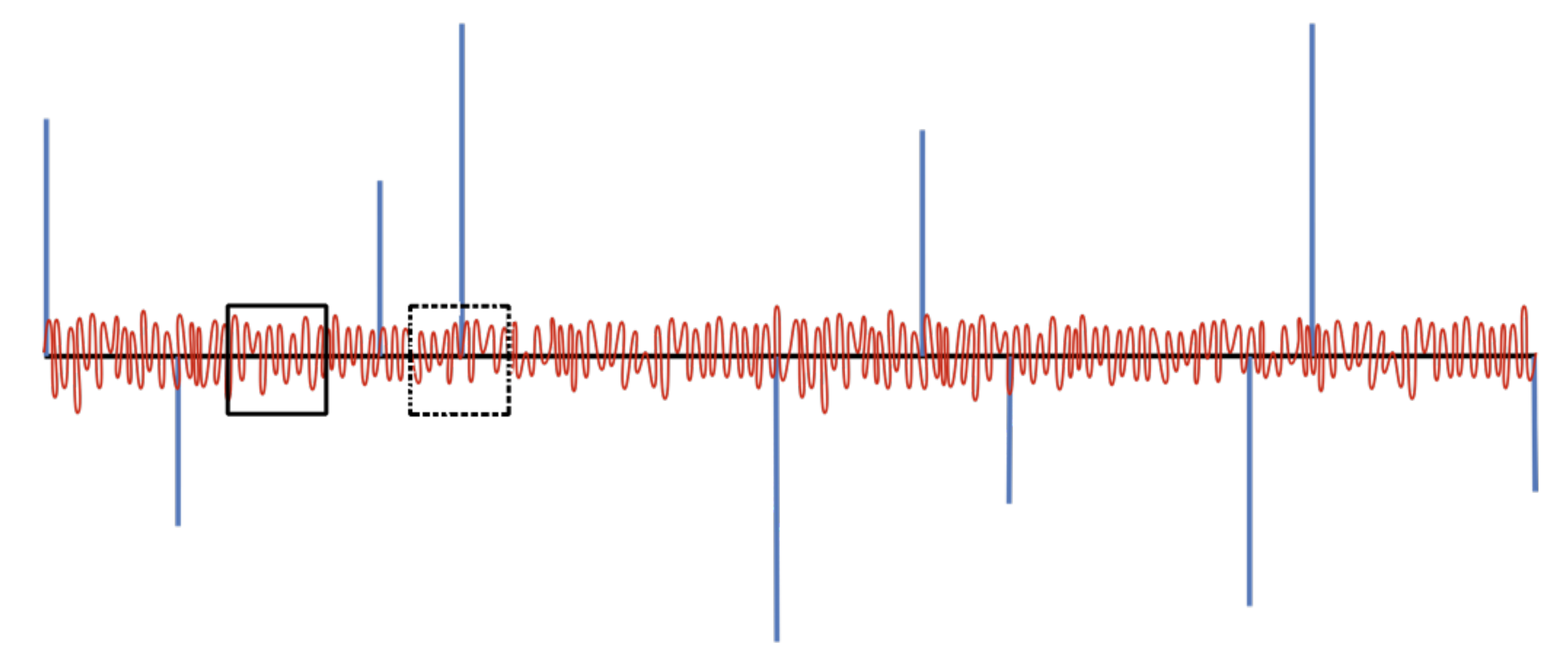}
\caption{Visualization of the Variance estimator in the setting of orthonormal design matrix. The received signal $y$ is a sum of a sparse signal (blue) and noise (red), and the variance estimate is an average of variances over ``good" windows (solid rectangle), but not ``bad" windows (dashed rectangle).}
\label{fig:spikes}
\end{figure}

\begin{algorithm}
\caption{Greedy Variance Estimator -- Orthonormal Design Matrix}
\label{GVE-1}

    \begin{algorithmic}[1]
    \State    Compute the window estimators $S_j = \frac{1}{L} \sum_{i\in \Omega_j} | y_i |^2, \quad \quad j \in \{1,2, \dots, p/L\}$.
    \State  Let $\widehat{\sigma^2} = \frac{2L}{p} \sum_{j=1}^{p/(2L)} S_{(j)}$, where $\{S_{(j)} \}_j$ is a non-decreasing arrangement of $\{S_j\}_j$.
    \end{algorithmic}

\end{algorithm}
\noindent The basic idea behind the above algorithm is that we want to capture a noise estimator that avoids the entries of $y$ affected by signal (hence in the second step we take the average of the smaller $50\%$ of the window estimates).   We can easily visualize the algorithm in this setting -- the received signal is a sum of a sparse signal and uniform noise, and the algorithm computes returns as an estimate of the noise variance the average of the smallst $50\%$ of the empirical variances taken over windows of length $L$.  We illustrate this in Figure \ref{fig:spikes}.

\noindent We have the following result which guarantees accuracy of the estimator $\widehat{\sigma^2}$.

\begin{theorem}
\label{Theorem_orth_mat}
Suppose $y = X \beta + \eta$ where $X \in \mathbb{R}^{p \times p}$ is orthonormal, $\eta_j \sim {\cal N}(0, \sigma^2)$ are independent, and $\beta$ is $s$-sparse.   Consider window size $L \geq \log^3(p)$, and suppose that $s \leq \frac{p}{2L}$.  Then the Greedy Variance Estimator produced by Algorithm \ref{GVE-1} satisfies
\begin{equation}
\label{grad_guarantee}
|   \widehat{\sigma^2} - \sigma^2 | \leq \frac{6}{\log p}\sigma^2,
\end{equation}
with probability $1 - \frac{2}{p}$.
\end{theorem}

\begin{remark}
(Total variation denoising) Suppose we receive image-type data and instead of taking the LASSO minimizer we want to instead want to regularize by the total variation seminorm:
\begin{equation}
\label{eq:tv_norm}
\widehat{\beta} = \text{arg}\min_\beta \| \beta - y\|_2^2 + 2 \lambda \text{TV}(\beta),
\end{equation}
where $\text{TV}(\beta):= \sum_n \| \beta_n - \beta_{n-1}\|$.  The TV penalty promotes sparsity with respect to the discrete derivative of signal, $\nabla \beta$.   In case $\nabla \beta$ is sparse, one may consider the natural analog of the greedy variance estimator of Algorithm \ref{GVE-1} where instead one takes window estimates of the form 
$$
S_j^{TV} = \frac{1}{L} \sum_{i\in \Omega_j} | y_{i+1} - y_i |^2, \quad \quad j \in \{1,2, \dots, p/L\}.
$$

Assuming that the discrete derivative vector $\nabla \beta$ is $s$-sparse, this modified greedy variance estimator provides a similar guarantee to \eqref{grad_guarantee} (up to a factor of 2) under the conditions of Theorem \ref{Theorem_orth_mat} on $X, L$, and $s$.  The only difference in the proof is that now, the random variable $y_{i+1} - y_i = (\beta_{i+1} - \beta_i) + (\eta_{i+1} - \eta_i)$ is not independent from its immediate neighbors $y_{i+2} - y_{i-1}$ and $y_i - y_{i-1}$, and so one applies the concentration results for i.i.d. variables separately to the even indices $2 i$ and the odd indices $2i + 1$, and combines the results with a standard union bound. 
\end{remark}
\bigskip

\section{Greedy Variance Estimation -- RIP Design Matrix}
\label{sect:rand_des}
We now turn to the more general case where the design matrix $X \in \R^{n\times p}$ is possibly underdetermined $n \leq p,$ but satisfies the Restricted Isometry Property with the appropriate constants (indeed this is a more general case, as an orthonormal matrix satisfies the RIP with constant $\delta = 0$).   We define the regularized design matrix as $Z := [Z_{\Omega_1},...,Z_{\Omega_{p/L}}]$ where each $Z_{\Omega_i} \in \R^{n\times L}$,
\begin{equation}
\label{Z_def}
    Z_{\Omega_i} := U_i I_{n\times L} V_i  \quad \text{such that} 
\end{equation}
\begin{equation*}
    X_{\Omega_i} = U_i \Sigma_i V_i \text{ is the SVD of }X_{\Omega_i}.
\end{equation*}
We first run a conditioning step based on the (block orthonormal) matrix $Z$ by computing $\tilde{y} = Z^T y$ before running the algorithm as in the orthonormal case with respect to the vector $\tilde{y}$.  The reason we apply $Z^T$ instead of the (faster to compute, and perhaps more natural) matrix $X^T$ is that $Z_{\Omega}^T X_{\Omega}$ will be better conditioned than $X_{\Omega}^T X_{\Omega}$ -- precisely, the singular values of the former will lie in the interval $[\sqrt{1-\delta}, \sqrt{1+\delta}]$ and the singular values of the latter in the interval $[1-\delta, 1 + \delta]$.   Applying $Z^T$ instead of $X^T$ not only improves the estimator quality in theory, but also in the numerical experiments.

\begin{algorithm}[H]
\caption{Greedy Variance Estimator}
\label{GVE-2}

\begin{algorithmic}[1]
    \State Compute $\tilde{y} = Z^T y$.
    \State Compute the window estimators $S_j = \frac{1}{L} \sum_{i\in \Omega_j} | \tilde{y}_i |^2,  \quad \quad j \in \{1,2, \dots, n/L\}.$
    \State Let $\widehat{\sigma^2} = \frac{2L}{p} \sum_{j=1}^{p/(2L)} S_{(j)}$, where $\{S_{(j)}\}_j$ is a non-increasing arrangement of the window estimators $\{S_j\}_j$.
\end{algorithmic}
\end{algorithm}

In practice, we use the matrix $X$ instead of $Z$, however using $Z$ allows us to do a more streamlined theoretical analysis.  To see why this should work intuitively, assume that we precondition just on $X$ that satisfies RIP for a large enough sparsity level $s_0$.  Note that $X^T y= X^T X \beta + X^T \eta$, so the obstruction to estimating the noise is the $X^T X$ term.  Then, $\|X \beta\|_2 = \| X_{\Omega_{\beta}} \beta \|_2 \approx \|\beta\|_2$, and if we assume our window set $\Omega_j$ is disjoint from $\Omega_{\beta}$, RIP implies the restricted matrices $X_{\Omega_j}^T$, $X_{\Omega_{\beta}}$ satisfy $\| X_{\Omega_j}^TX_{\Omega_{\beta}}\| \leq \delta$ for $\delta>0$ small.  Thus, for a ``good" window estimator, we only see the noise $X^T \eta$.  The constants in Theorem \ref{Rand_mat_est} are chosen for neatness of presentation and are in no way optimized.

\begin{theorem}
\label{Rand_mat_est}
Suppose $y = X \beta + \eta$ where $X \in \mathbb{R}^{n \times p}$  has the RIP of order $s$ and level $\delta$, $\eta_j \sim {\cal N}(0, \sigma^2)$ are i.i.d., and $\beta$ is $s$-sparse.   Assume that $\log^3(p) \leq L \leq n$, and $2L \leq s \leq \frac{n}{2L}$.   Then, the Greedy Variance Estimator in Algorithm \ref{GVE-2} satisfies
    \begin{equation*}
    \left| \widehat{\sigma}^2 - \sigma^2 \right|  \leq   \frac{2 \delta \| \beta \|^2}{L} + \frac{5 \sigma^2}{\log(p)} + \frac{1}{L} \max\left( 4\sigma^2 \log(p), 8 \sqrt{\delta}\sigma \| \beta \|^2 \sqrt{\log(p)}  \right) 
    \end{equation*}
    with probability $1- \frac{4}{p}$.
\end{theorem}
To make this result more concrete, we provide a corollary of Theorem \ref{Rand_mat_est} for the case that  $X \in \mathbb{R}^{n \times p}$ is a Gaussian matrix with i.i.d. entries $X_{j,k} \sim {\cal N}(0,1/n)$.   Recall \cite{foucart2013mathematical, baraniuk2008simple} that for an $n \times p$ matrix $X$ with i.i.d. entries $X_{i,j} \sim {\cal N}(0,1/n),$ the RIP of order $s$ and level $0 < \delta < 1$ is satisfied with probability at least $1- 2\exp(-\delta^2 n/(2C))$ once
\begin{equation}
 C \frac{s\log(p/s) }{n} \leq \delta^2 < 1
\label{s}
\end{equation}
for a universal constant $C > 0$.  In the following corollary, we set $\delta^2 =  C \frac{s\log(p/s) }{n} < 1$ and $L = \frac{n}{2s}$.

\begin{corollary}
\label{cor:mainThm}
Fix integers $s < n < p$ such that $Cs \log^3(p)/n < 1 \leq s^2/n$.  Draw $X \in \mathbb{R}^{n \times p}$ as a Gaussian matrix with i.i.d. entries $X_{j,k} \sim {\cal N}(0,1/n)$.  
Suppose $y = X \beta + \eta$ where $\eta_j \sim {\cal N}(0, \sigma^2)$ are i.i.d., and $\beta$ is $s$-sparse.   Then, the Greedy Variance Estimator in Algorithm \ref{GVE-2} with window size $L=\frac{n}{2s}$ satisfies
    \begin{equation*}
    \left| \widehat{\sigma}^2 - \sigma^2 \right|  \leq   \frac{2 s \| \beta \|^2}{n} + \frac{5 \sigma^2}{\log(p)} + \frac{8 \sigma^2 s \log(p)}{n} + \frac{16 s \sigma \| \beta \|^2 \sqrt{\log(p)} }{n}
    \end{equation*}
    with probability at least $1- \frac{4}{p} - 2 (s/p)^s$, with respect to both the draw of $X$ and the draw of $\eta$.  Here, $C > 0$ is a universal constant independent of all other parameters.
\end{corollary}

\begin{remark}
{\bf Consistency.}  Our estimator is asymptotically consistent under certain sequences $\beta = \beta_p \in \mathbb{R}^p$.  For instance, fix $\rho < 1$, and consider values of $n \leq p$ such that $p \leq \rho^{-1} n$.  Consider the assumptions of Corollary \ref{cor:mainThm}, and moreover set $s = \sqrt{p}$.  Then, for any sequence $\beta = \beta_{p}$ such that $\| \beta_{p} \|_0 \leq s$ and $| \beta{p_{(j)}} | \leq \gamma $ for each $j \leq s$, we have with probability going to 1,
\begin{equation}
\label{error:asymp}
 \sup_{0 \leq \sigma^2 \leq c} \mathbb{E}_{\eta}  | \widehat{\sigma}^2 - \sigma^2 | \leq C\frac{\gamma^2 c}{\rho^2 \log(p)}
\end{equation}
In particular, this implies that the greedy variance estimator is consistent in the sense that restricted to such sequences $\beta_p$, 
$$
\lim_{n/p \rightarrow \rho} \sup_{0 \leq \sigma^2 \leq c} \mathbb{E} | \sigma^2 - \hat{\sigma}^2 | = 0.
$$
Note that the estimator achieves asymptotic consistency even as $\| \beta \|_1$ is allowed to grow as $\gamma \sqrt{p}$.  
By comparison,  Theorem 2.1 of \cite{chatterjee2015prediction} concerning asymptotic consistency of a variant of cross-validated LASSO requires that $\|\beta \|_1 \leq \frac{n^{1/4}}{\log(p)^{1/4}}$, but does not require that $\beta$ is $\sqrt{p}$-sparse.  Alternatively, Lemma 2 in \cite{dicker2014variance} proves consistency for the method of moments estimator in case the design matrix is Gaussian, and under the stronger assumption that $\lim_{p \rightarrow \infty} \| \beta_p \|_2 \leq c$ remains bounded.  
\end{remark}

\section{Experiments}
\label{sect:experiments}
Our experimental methodology is based off of the results in \cite{reid2013study}.  In particular, we generate a design matrix $X\in \R^{n\times p}$ with i.i.d. entries $X_{ij} \sim \mathcal{N}(0,n^{-1/2})$ so that $X$ satisfies RIP with sufficiently small constants with high probability.  The sparsity level $s = \lceil n^{\alpha} \rceil$, with $\alpha < 1$, and the non-zero entries of $\beta$ (chosen uniformly at random) are distributed according to a Laplace(1) distribution.  The resulting $\beta$ is scaled to have the specified norm.
The experiments are over the following grid of parameter values, where $n=100$ in all experiments.
\begin{itemize}
    \item $p= 100,200,500,1000,$
    \item $\|\beta\|_2 = 0.1,1,2,5,10,$
    \item $\alpha = 0.1,0.3,0.5,0.7,0.9.$
\end{itemize}

We use the following estimators in our analysis:
\begin{itemize}
    \item oracle: the oracle estimator $\hat{\beta} = \|\eta\|_2/\sqrt{n}$.
    \item window: the standard window estimator with the transformation $\tilde{y} = X^Ty$.
    \item window-svd:  the theoretical window estimator with the transformation $\tilde{y} = Z^T y$ where $Z$ is given by \eqref{Z_def}.
   
    \item ${ \text{window} }_{BC}$ and ${ \text{window-svd} }_{BC}$: same as window and window-svd, except we multiply the resulting estimators by $1+ \frac{1}{\log(p)}$ to correct the downward bias that results from averaging only over the windows with the smallest statistics.  This estimator could be preferable in practical applications where it is crucially important that the estimator does not \emph{underestimate} the variance. 
    
    \item cv-lasso: 10-fold cross-validated LASSO (computed using the R package glmnet \cite{glmnet}).
    \item moment: method of moments estimator (see \cite{dicker2014variance}).
\end{itemize}
We include the cross-validated LASSO because it was shown to be the most robust to changes in sparsity/dimension by \cite{reid2013study} and the method of moments estimator because it aims to be a fast replacement for cv-LASSO.
The window size is chosen based on an inflection point in the values of the estimator for a specific set of parameters as the window size varies.\\

As we can see in Figure \ref{fig:oblivious}, the window and window-svd estimators have reasonable performance compared to the cv-LASSO with slightly larger biases.  In particular, we do quite well for $\alpha=0.1$, $\beta=1$, performing similarly to cv-Lasso, and with a much smaller variance than the method of moments.

\begin{remark}
We only include results for $\alpha = 0.1$ because the algorithm performs similarly for $\alpha \leq 0.5$.  Moreover our theory only covers up to roughly $\alpha = 0.5$ for reasonable choices of window size.  The performance for dense signal $\alpha = 0.9$ is covered in its own section below.
\end{remark}

\subsection{Optimal Window Size}
\label{subsect:opt_window_size}
It is notable to see how well our method can perform when the window size is optimized.  Here, we give some representative plots (Figures \ref{fig:opt} and \ref{fig:opt_10000}) to show what happens to performance when replacing the window size with the optimal window size using prior knowledge of the variance.  In all experiments, n=100 and p=1000,10000.  For the low SNR regimes, we see a similar downward bias to the oblivious choice of window size, although with a smaller bias.  Similarly, for high SNR, the upward bias is also smaller than when choosing an oblivious window size.  By scaling our plots to p=10000, we see consistent results in higher dimensions.  In table \ref{table:opt_wind_size} we report the optimal window size for various values of $\alpha$ and $\|\beta\|_2$.  The optimal window size was found by a grid search over all possible window sizes using knowledge of the true variance.

\begin{table}[h!]

\centering
\begin{tabular}{|c|c||c|c|c|c|c|}
     \cline{3-7}
     \multicolumn{2}{c||}{}&\multicolumn{5}{|c|}{$\|\beta\|_2$}\\
     \cline{3-7}
     \multicolumn{2}{c||}{}&0.1&1&2&5&10\\
     \cline{3-7} \hline \hline
     \multirow{4}{.25cm}{$\alpha$} &0.1& 100 & 100 & 100 & 20&4\\
     \cline{2-7}
     &0.3&100&100&100&22&4\\
     \cline{2-7}
     &0.5&100&100&100&18&4\\
     \cline{2-7}
     &0.7&100&100&100&18&4\\
     \cline{2-7}
     &0.9&100&100&100&14&3\\
     \hline
\end{tabular}
\caption{Optimal window sizes as a function of $\alpha$ and $\|\beta\|_2$ for $p=200$.  We note that the optimal window size is generally decreasing as a function of both the signal to noise ratio and the sparsity.  Moreover, choosing the maximal window size is optimal in modest regimes.}
\label{table:opt_wind_size}
\end{table}

\subsection{High Dimensional Regime}
In this section we highlight the regime in which our estimator is most useful - when $p \gg n$ is large.  In particular, we chose $n=100$, $p=100000$ in all experiments.  In this regime, it is inefficient to even compute an optimal box size based on an inflection point in the value of the estimator, so instead the choice $L=25$ was fixed for all experiments.  The results are shown in Figure \ref{fig:high}.  In low SNR regimes we see the bias corrected estimators (window(BC) and window-svd(BC)) outperforming the normal estimators.  Although the bias remains, the estimator performs well, especially in low SNR regimes.  This is likely due to the strength of the compressed sensing properties for the design matrix as the dimension grows.  The bias increases with higher SNR, however our estimator maintains a lower variance than cv-LASSO.

\subsection{Orthogonal Design Matrix}

We find our estimator performs quite well in the case where the design matrix is orthogonal, as shown in Figure \ref{fig:ortho}.  In all experiments, $p=n=200$ and the window size is chosen via inflection point in the value of the estimator. The method of moments still performs reasonab   ly well, but suffers a strong upwards bias for large SNR.  We note that in all regimes, our estimator performs better than cross-validated LASSO.  Moreover, it is more robust to changes in SNR than when the design matrix is RIP (but not necessarily orthogonal).

\subsection{Dense Signal}
Our theory does not cover high sparsity levels ($\alpha \geq 0.9$), but nonetheless our estimator performs well.  Although more prone to high levels of SNR, we are still competitive with cv-LASSO in low SNR regimes as seen in Figure \ref{fig:dense}.

 \subsection{Real Data}
In this section we report results on real data sets well suited for LASSO, but whose design matrices are far from satisfying the RIP.  We illustrate that nevertheless, our greedy variance estimator works comparably to CV-LASSO for such data sets.  Typical data sets where $p \gg n$ involve genetics data, where the amount of genetic data recorded is much larger than the number of patients sampled.

The first data set is from \cite{singh2002gene} and corresponds to gene expression data.  It is presented as a $102 \times 6033$ matrix, where each row is a sample from a single subject, and the columns are expression levels.  We defer to the original paper for how precisely these values were computed.  This data is regressed against a length 102 vector with 52 cancer patients (1) and 50 healthy patients (0).  We also consider the well-known Golub data set \cite{golub1999molecular}, which is a gene expression data set from subjects with human acute myeloid (AML) and
acute lymphoblastic leukemias (ALL).  It is represented as 3571 expression levels over 72 patients, with 47 ALL subjects and 25 AML.  The final data set is from Alon et al. \cite{alon1999broad}, a 62x2000 matrix of gene expression data from colon tissue, 40 tumor 22 normal.  Note that in all cases we have a small number of subjects ($<102$) and thousands of gene expressions for each subject.

Since for real world data we cannot compare our estimated variance to the ``true" underlying noise level, we instead compare the noise variance computed for 10 fold CV-LASSO to that returned by the greedy variance estimators, as well as the resulting $\lambda$ parameters.  These results are tabulated in table \ref{table:lasso}.  Here,
\begin{itemize}
\item $\sigma$ / $\lambda$ GVE corresponds to the window-SVD estimates of $\sigma$ and $\lambda$, respectively
\item $\sigma$ / $\lambda$ Fast-GVE corresponds to the window estimates of $\sigma$ and $\lambda$, respectively
\item $\lambda$ min-MSE corresponds to the value of $\lambda$ obtained via CV-LASSO corresponding to minimal MSE on the hold-out data
\item $\lambda$ 1-SE corresponds to the value of $\lambda$ obtained via CV-LASSO which is one standard error below $\lambda$ min-MSE, as often used to avoid overfitting
\end{itemize}
In the experiments in this section, we chose window size $L = 2 \log(p)$ which was not optimized in any way.

In table \ref{table:lasso}, we note that for GVE, both the estimated variance and resulting $\lambda$ parameter are close to the corresponding $\lambda$ value for 1 standard error in CV-LASSO.

We also plot, in figures \ref{fig:singh}-\ref{fig:leukemia} the corresponding curves for the mean squared error of the LASSO solution, using the $\lambda$ parameters from table \ref{table:lasso}.

\begin{table}
    
    \centering
    \begin{tabular}{|c|c|c|c|c|c|c|}
        \hline
         Data& $\sigma$ CV-LASSO & $\sigma$ GVE & $\sigma$ Fast-GVE & $\lambda$ 1-SE & $\lambda$ min-MSE & $\lambda$ GVE \\
         \hline
         \cite{singh2002gene}&0.4854	&0.7254	&105.279	&0.05295	& 0.0954&	0.00429\\
         \hline
          \cite{golub1999molecular}& 0.8132   &	0.6772 & 375.82	&0.0637	&0.0473& 0.0276\\
          \hline
          \cite{alon1999broad}&0.7788 & 1.212 &8.31E+09	& 0.1503	& 0.1699	&0.0861\\
          \hline
    \end{tabular}
    \caption{$\sigma$ and $\lambda$ values for the genomics data sets.}
    \label{table:lasso}
\end{table}

\section{Future Work}
\label{sect:future}
Our estimator has been shown to be a useful tool to use in high dimensional variance estimation, and comes with nice theoretical properties that leverage results from the compressed sensing literature.  Moreover, it is extremely fast, parallelizable, and is competitive with cv-LASSO in most parameter regimes. Based on our experimental/theoretical results there are some obvious directions to go in the future:
\begin{itemize}
    \item Develop an efficient estimator that has theoretical guarantees for a more general design matrix, in particular that satisfies the restricted eigenvalue condition.
    
    \item Find a choice of box size that is more robust to sparsity and SNR, which is still efficient to compute.
\end{itemize}
Although this estimator is by no means a replacement for existing estimators in typical regimes, it scales extremely well into high dimensions and performs as well if not better when $p \gg n$.  This regime seems the most interesting for developing more robust estimators.

\section*{Acknowledgments}
We thank Abhinav Nellore for discussions on parameter selection in high dimensional problems which motivated this work.  We also thank Robert Tibshirani for directing us to the glmnet package for computing cv-LASSO.  We also thank the anonymous referees for their feedback which greatly improved the manuscript.  R. Ward and C. Kennedy were partially supported during this work by NSF CAREER grant \#1255631.

\bibliographystyle{alpha}
\bibliography{GRERevision}

\appendix

\section{Proof Ingredients}
\begin{proposition} (Lemma 1 in \cite{laurent2000adaptive}) Suppose $Z$ has a chi-squared distribution with $d$ degrees of freedom.  Then,
\begin{equation}
    \label{Chi_Square_Conc}
    \P[d - 2 \sqrt{dt} \leq Z \leq d + 2\sqrt{dt} + 2t] \geq 1 - 2 e^{-t}\quad \quad \quad \forall t \geq 0.
\end{equation}
\end{proposition}

\begin{proposition}(Proposition 2.5 in \cite{rauhut2010compressive}) Suppose $\Omega_u \cap \Omega_v = \varnothing$, and that $X \in \R^{n \times p}$ satisfies RIP of order $s_0$ and level $\delta >0$ with $s_0 = |\Omega_u| +  |\Omega_v|$.  Then,
\begin{equation}
\label{DisjointRip} \| X_{\Omega_u}^T X_{\Omega_v} \|_{2 \rightarrow 2} \leq \sqrt{\delta}
\end{equation}

\end{proposition}

\begin{proposition}
\label{prop:gaussianconc}(Equation (5.5) in \cite{vershynin2010introduction})
Let $X$ be a Gaussian random variable with mean 0, variance $\sigma$.  Then, 
\[
\P[|X|>t] \leq 2 e^{-t^2/2\sigma^2}, \hspace{1cm} t\geq 1.
\]
\end{proposition}

\section{Proofs}
\subsection{Proof of Theorem \ref{Theorem_orth_mat}}

Consider the window estimators 
\begin{align}
S_j &= \frac{1}{L} \sum_{i\in \Omega_j} | y_i |^2 \nonumber \\
&= \frac{1}{L}  \sum_{i\in \Omega_j} | \beta_i + \eta_i |^2 \nonumber \\
&= \frac{1}{L}  \sum_{i \in \Omega_j} | \beta_i |^2   + \frac{1}{L}  \sum_{i \in \Omega_j}   |  \eta_i |^2 + 2 \frac{1}{L}  \sum_{i\in \Omega_j} \beta_i \eta_i.  \nonumber
\end{align}
Set $E_j := \frac{1}{L}  \sum_{i \in \Omega_j}  |  \eta_i |^2$. 
$E_j$ is a sum of $L$ independent squares of ${\cal N}(0,\sigma^2)$ random variables.  Then $E_j$ concentrates strongly around its expected value, 
$$\mathbb{E}(E_j) = \sigma^2.$$
Note that $E_j$ has a chi-squared distribution with $L$ degrees of freedom, so by \eqref{Chi_Square_Conc}
with the choice $t = \log(p^2) $ and after a union bound over all $p/L$ windows, we get that with probability at least $1-\frac{2}{p}$,
$$
\left(1 - \frac{5}{\log(p)}\right) \sigma^2 \leq E_j \leq \left(1 + \frac{5}{\log(p)}\right) \sigma^2, \nonumber
$$
holds uniformly for all $j \in \{1,2, \dots, p/L\}$, assuming that $L \geq \log^3(p)$.  

  Since $L \leq \frac{p}{2s}$ by assumption, the pigeon hole principle implies that at least $\frac{p}{2L}$ windows do not overlap $\Omega_{\beta}.$
  On any such ``good" window $k$ we have $\| \beta_{Lk:Lk+L-1} \|_2^2 = 0$ and hence
  \begin{align}
| S_{k} - \sigma^2|  &\leq  \frac{5\sigma^2}{\log(p)} .
  \end{align}
  Thus, if $\overline{S}$ is the average over a subset of the good windows, then also $| \overline{S} - \sigma^2 | \leq \frac{5 \sigma^2}{\log(p)}$. 
  
Now, to bound the estimator above on \emph{any} window, we need some control on the cross term $ \sum_{i\in \Omega_j} \beta_i \eta_i$.  Note that this quantity is just a sum of i.i.d. Gaussians with mean zero and with variance $ \| \beta_{\Omega_j \cap \Omega_{\beta}}\|_2^2 \sigma^2$; thus, by concentration, we have that with probability at least $1 - 2/p$, the following holds uniformly over all windows:
\begin{equation}
\frac{1}{L} \sum_{i\in \Omega_j} \beta_i \eta_i \leq \frac{2 \sigma \|\beta_{\Omega_j \cap \Omega_{\beta}} \|_2 \sqrt{\log (p)}}{L}.
\end{equation}
Hence, for any \emph{any} window,
  \begin{align}
 S_{j}  &\geq  \frac{1}{L} \| \beta_{\Omega_j \cap \Omega_{\beta}} \|_2^2  +  E_j -  \frac{2}{L} \sum_{i \in \Omega_j \cap \Omega_{\beta}}  \beta_i \eta_i  \nonumber \\
&\geq  \frac{1}{L} \| \beta_{\Omega_j \cap \Omega_{\beta}} \|_2^2  + \left(1-\frac{5}{\log(p)}\right) \sigma^2  -  \frac{\| \beta_{\Omega_j \cap \Omega_{\beta}} \|_2}{\sqrt{L}} \frac{2\sigma \sqrt{\log{p}}}{\sqrt{L}}    \nonumber \\
    &\geq
    \left(1-\frac{5}{\log(p)}\right) \sigma^2 -  \frac{ \sigma^2 \log(p)}{L}  \quad \quad \quad \text{The minimal value of a quadratic } x^2 + b - ax  \text{ is } b - a^2/4 \nonumber \\
        &\geq  \left(1-\frac{6}{\log(p)}\right) \sigma^2 \nonumber,
         \end{align}
where the final inequality holds because $\log^2(p)) \leq L$.

Now, consider the estimator $\widehat{\sigma^2} = \frac{2L}{p} \sum_{j=1}^{p/(2L)} S_{(j)}$.  By construction, $\widehat{\sigma^2} \leq \overline{S}$, where $\overline{S}$ is the average over any $p/(2L)$ ``good" windows.  From the above analysis, we have that with probability exceeding $1- \frac{4}{p}$,
 $$
 | \widehat{\sigma^2} - \sigma^2| \leq \frac{6}{\log(p)} \sigma^2.
 $$
 
\subsection{Proof of Theorem \ref{Rand_mat_est}}
Recall that $\tilde{y}:= Z^T y \in \mathbb{R}^p$.  Consider the window estimate
\begin{align}
    S_j &= \frac{1}{L} \sum_{i \in \Omega_j} |\tilde{y}_i|^2 \nonumber\\
    &=  \frac{1}{L} \sum_{i \in \Omega_j} |(Z^T X \beta)_i|^2 + \frac{1}{L} \sum_{i\in \Omega_j} |Z_i^T \eta|^2 + \frac{2}{L} \sum_{i \in \Omega_j}  (Z^T X \beta)_i (Z^T \eta)_i\nonumber\\
    &= \frac{1}{L} \| Z_{\Omega_j}^T X_{\Omega_{\beta}} \beta\|_2^2 + \frac{1}{L} \| Z_{\Omega_j}^T \eta \|_2^2 + \frac{2}{L} \sum_{i \in \Omega_j}  (Z^T X \beta)_i (Z^T \eta)_i \label{window_square_exp}
\end{align}
The first term is small if $\Omega_j$ and $\Omega_{\beta}$ have disjoint support, since $X$ has the RIP of order $s \geq 2L$.   The center term gets close to its expectation $\sigma^2$ due to standard concentration inequalities, and the third term is also small due to standard concentration inequalities.  More concretely, if we assume that $S_j$ is a ``good" window, meaning that $\Omega_j$ and $\Omega_{\beta}$ have disjoint support, by equation \eqref{DisjointRip}
\begin{equation}
\label{first_term_x}
    \frac{1}{L} \| X^T_{\Omega_j} X_{\Omega_{\beta}} \beta \|_2^2 \leq \frac{\delta \|\beta\|_2^2}{L}.
\end{equation}
All of the diagonal entries of $\Sigma_j$ are in the range $[\sqrt{1-\delta},\sqrt{1+\delta}]$, hence by \eqref{first_term_x}
\begin{align}
\frac{1}{L} \| Z^T_{\Omega_j} X_{\Omega_{\beta}} \beta \|_2^2 &\leq \frac{1 + \delta}{L} \| X^T_{\Omega_j} X_{\Omega_{\beta}} \beta \|_2^2 \nonumber\\
&\leq \frac{\delta(1 + \delta) \| \beta\|_2^2}{L} \nonumber\\
&\leq \frac{2 \delta \|\beta\|_2^2}{L} \label{first_term}
\end{align}

For the center term, note that due to rotation invariance of the Gaussian, $\|Z_{\Omega_j} \eta\|_2^2 = \|P_L \tilde{\eta} \|_2^2$  where $\tilde{\eta}$ is an i.i.d. Gaussian vector, and $P_L$ is a projection onto the first $L$ coordinates. 
Next, we know that $\| P_L \tilde{\eta} \|_2^2$ has a chi-squared distribution with $L$ degrees of freedom, so by \eqref{Chi_Square_Conc} with $t = \log(p^2)$,
$$
\P\left[ | \| P_L  \tilde{\eta} \|_2^2 - L \sigma^2 | \leq 2 \sigma^2\left(\sqrt{L \log (p)} +  \log(p)\right) \right]\geq 1 - \frac{2}{p^2}.
$$
Hence by a union bound, with probability at least $1- \frac{2}{p}$, the following holds uniformly over all windows:
\begin{align}
    | \| Z^T_{\Omega_j} \eta \|_2^2/L - \sigma^2| & = | \| P_L  \tilde{\eta}  \|_2^2 /L -  \sigma^2 |\nonumber \\
    &\leq  2 \sigma^2\sqrt{\frac{\log(p)}{L}} + 2 \sigma^2 \frac{\log(p)}{L} \nonumber\\
    & \leq  \frac{5 \sigma^2}{\log(p)} \label{window_conc_ineq}
\end{align}

For the final term in \eqref{window_square_exp}, note that $\frac{2}{L} \sum_{i \in \Omega_j} (Z^T X \beta)_i (Z^T \eta)_i$ is a Gaussian random variable with variance $2\sigma \|Z^{T}_{\Omega_j} X_{\Omega_\beta} \beta \|_2 / L$.  Thus, by Proposition \ref{prop:gaussianconc} and \eqref{first_term}, the following holds uniformly over all windows with probability at least $1 - \frac{1}{p}$:
\begin{align}
    \frac{2}{L} \sum_{i \in \Omega_j} (Z^T X \beta)_i (Z^T \eta)_i & \leq \frac{4 \sigma \| Z_{\Omega_j} X_{\beta} \beta\|_2 \sqrt{\log(p)}}{L} \label{third_term1} \\
    &\leq \frac{8  \sqrt{\delta}\sigma \| \beta\|_2 \sqrt{\log(p)}}{L} \label{third_term},
\end{align}

Thus, averaging over any set of $p/2L$ ``good" windows, using \eqref{first_term} \eqref{window_conc_ineq} and \eqref{third_term} we have
\begin{align} 
   \left| \frac{2L}{p} \sum_j S_j  - \sigma^2 \right| &\leq \frac{2 \delta \|\beta\|_2^2}{L}+ \frac{5\sigma^2}{\log p} +  \frac{8\sqrt{\delta} \sigma \| \beta\|_2 \sqrt{\log p}}{L}
   \label{good_window_ineq}
 \end{align}
 with probability at least $1 - \frac{4}{p}$.   Thus, by construction, the estimator $\widehat{\sigma^2} = \frac{2L}{p} \sum_j S_{(j)}$ also satisfies 
 $$
\widehat{\sigma^2} \leq \sigma^2 + \frac{2 \delta \|\beta\|_2^2}{L}+ \frac{5\sigma^2}{\log p} +  \frac{8\sqrt{\delta} \sigma \| \beta\|_2 \sqrt{\log p}}{L}
 $$

It remains to show that the window estimator $\widehat{\sigma^2}$ cannot be too small.   The inequalities \eqref{third_term1} and \eqref{window_conc_ineq} hold uniformly over all windows, not just good windows; hence, for any window $S_j$,
\begin{align*}
    S_j &\geq \frac{1}{L} \| Z^T_{\Omega_j} X_{\Omega_{\beta}} \beta \|_2^2 + \frac{1}{L} \| Z_{\Omega_j}^T \eta \|_2^2  - \frac{2}{L} \| Z^T_{\Omega_j} X_{\Omega_{\beta}} \beta \|_2 \| X_{\Omega_j}^T \eta\|_2 \\
    &\geq \frac{1}{L} \| X_{\Omega_j}^T X_{\Omega_{\beta}} \beta \|_2^2 +\sigma^2 -  \frac{5\sigma^2}{\log(p)}   - \frac{8\sigma  \| X^T_{\Omega_j} X_{\Omega_{\beta}} \beta \|_2 \sqrt{\log p}}{L}\\
    &\geq \sigma^2 -\frac{5\sigma^2}{\log(p)}  - \frac{4 \sigma^2 \log(p)}{L}   \quad \quad \quad \text{(the minimal value of a quadratic } x^2 + b - ax  \text{ is } b - a^2/4)
\end{align*}
Combining the bounds,
$$
-\frac{5\sigma^2}{\log(p)}  - \frac{4 \sigma^2 \log(p)}{L}   \leq \frac{2L}{p} \sum_j S_{(j)} - \sigma^2  \leq  \frac{2 \delta \|\beta\|_2^2}{L}+ \frac{5\sigma^2}{\log p} +  \frac{8\sqrt{\delta} \sigma \| \beta\|_2 \sqrt{\log p}}{L}
$$
Thus, 
$$
| \widehat{\sigma^2} - \sigma^2 | \leq \left( 2\delta \frac{\| \beta \|^2}{L} + \frac{5 \sigma^2}{\log(p)} + \frac{1}{L} \max\left( 4\sigma^2 \log(p), 8 \sqrt{\delta}\sigma \| \beta \|^2 \sqrt{\log(p)}  \right)  \right)
$$

\pagebreak

\begin{figure}
    \centering
  \includegraphics[width=0.48\linewidth]{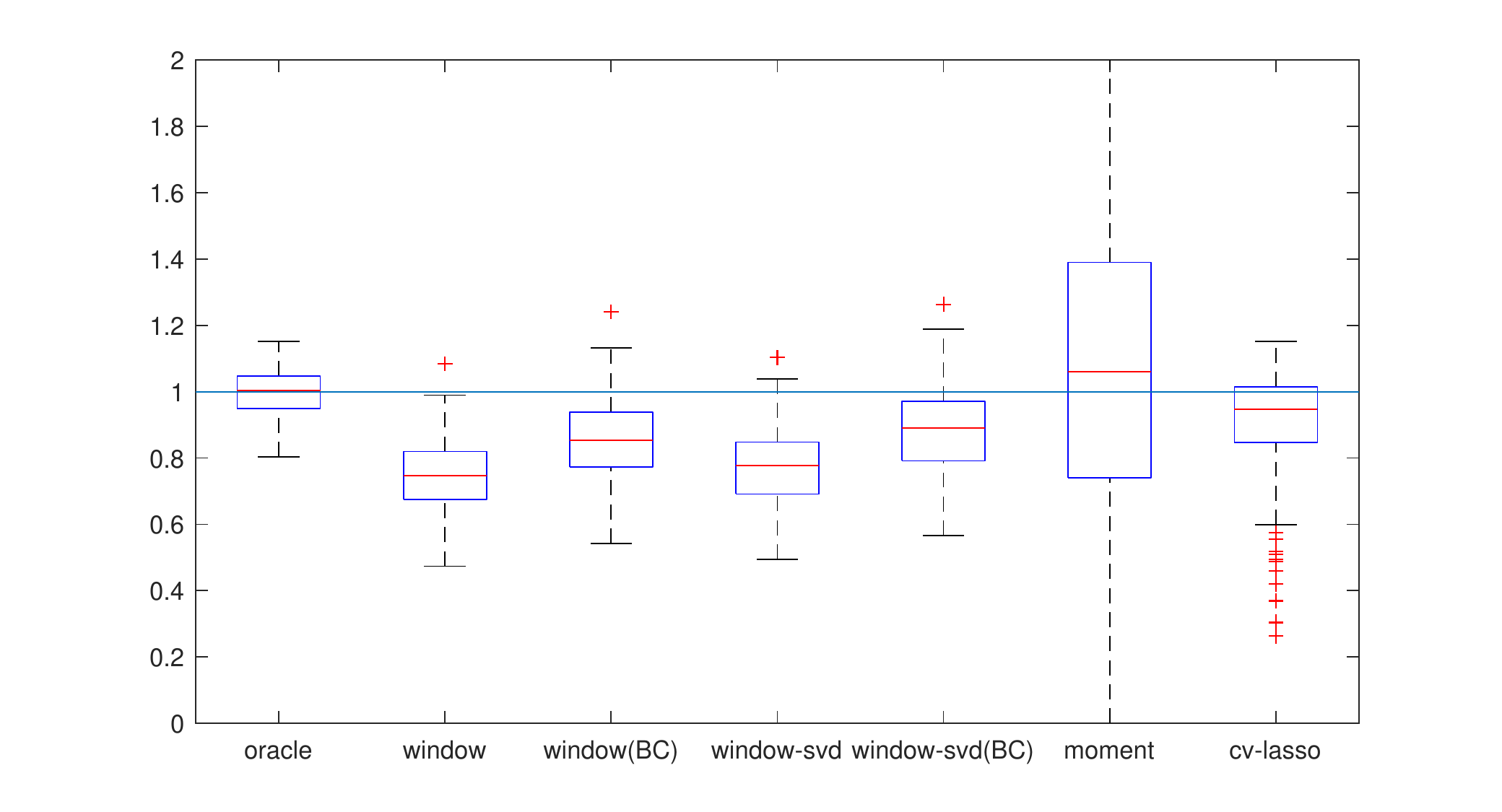}
  \includegraphics[width=0.48\linewidth]{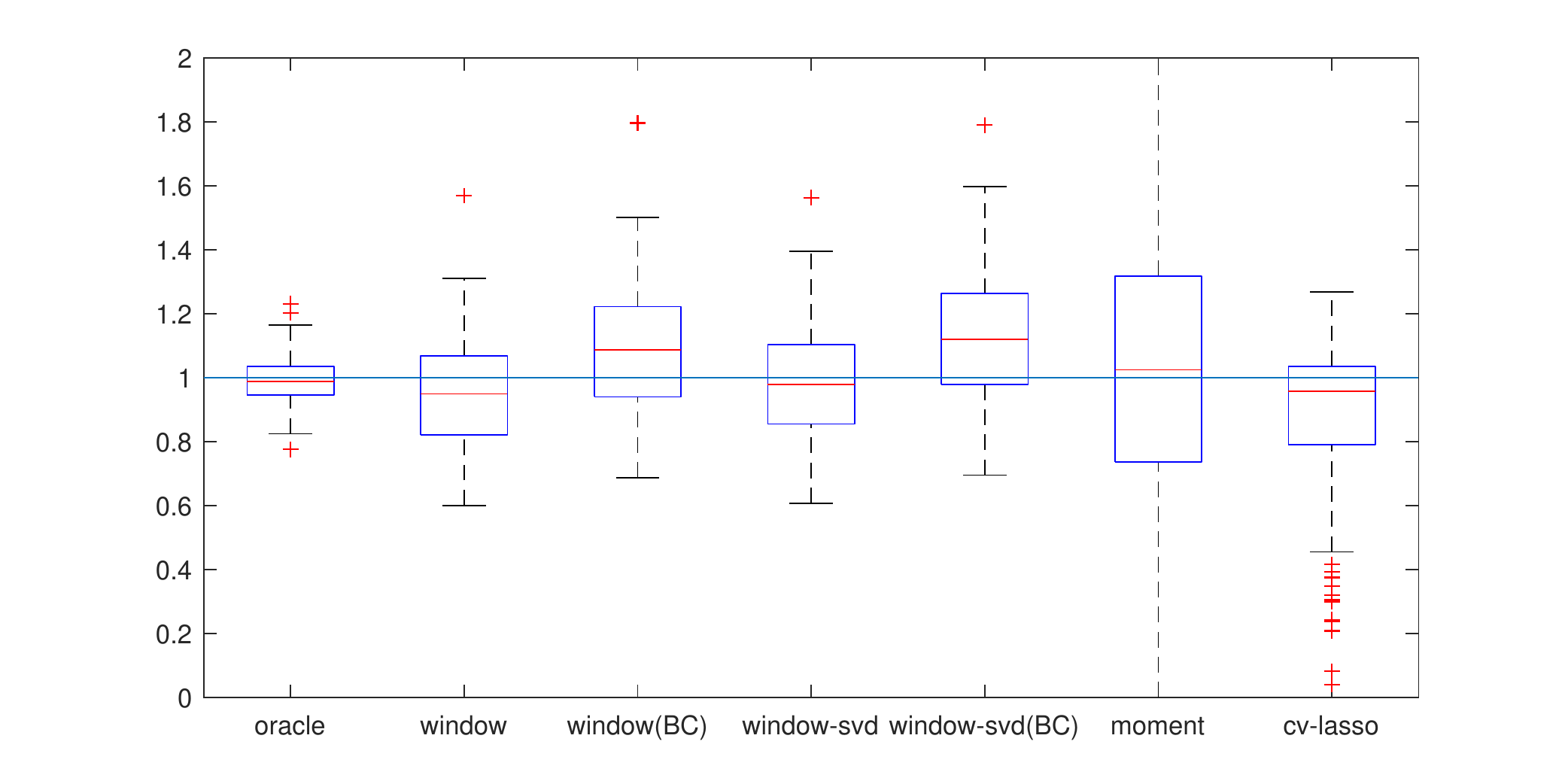}\\ 
    \includegraphics[width=0.48\linewidth]{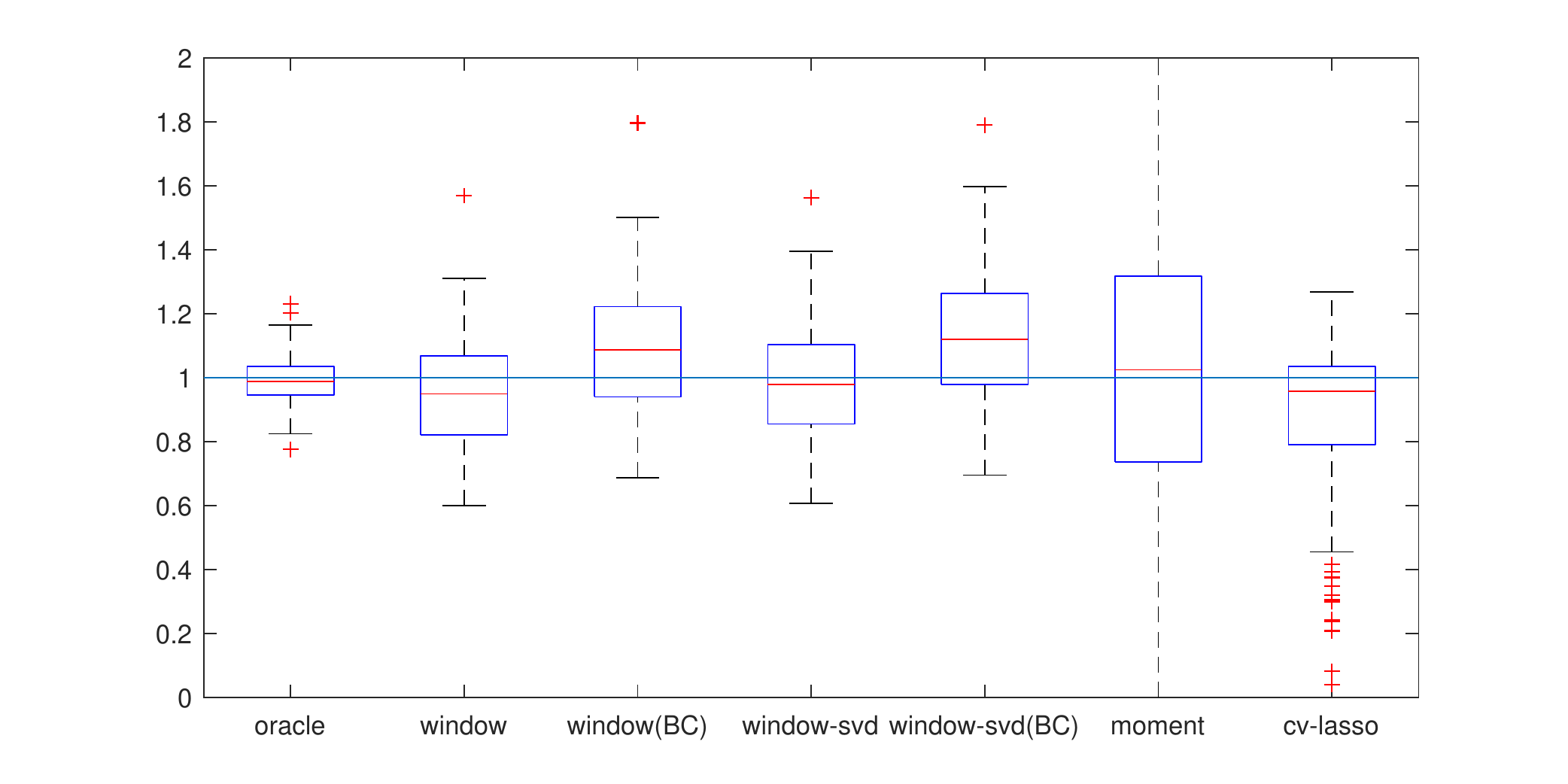}
     \includegraphics[width=0.48\linewidth]{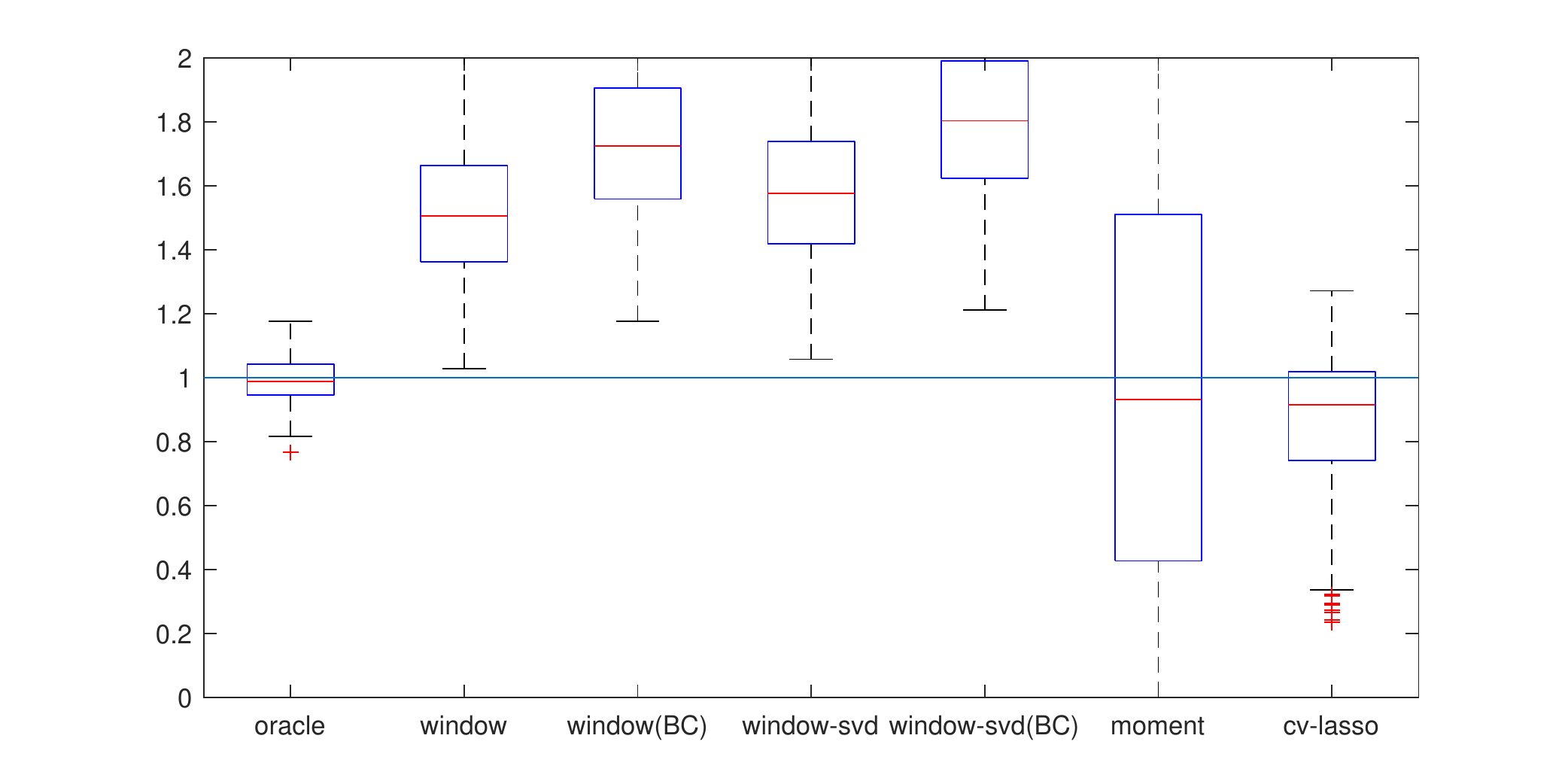}   
    \caption{Window size chosen based on inflection point, $p=1000$.  Signal-less ($\|\beta\|=0$), low SNR ($\alpha = 0.1$, $\|\beta\| = 1$), medium SNR ($\alpha = 0.1$, $\|\beta\| = 5$), high SNR ($\alpha = 0.1$, $\|\beta\| = 10$) respectively, top to bottom.}
    \label{fig:oblivious}
\end{figure}

\begin{figure}
    \centering

\includegraphics[scale=.55]{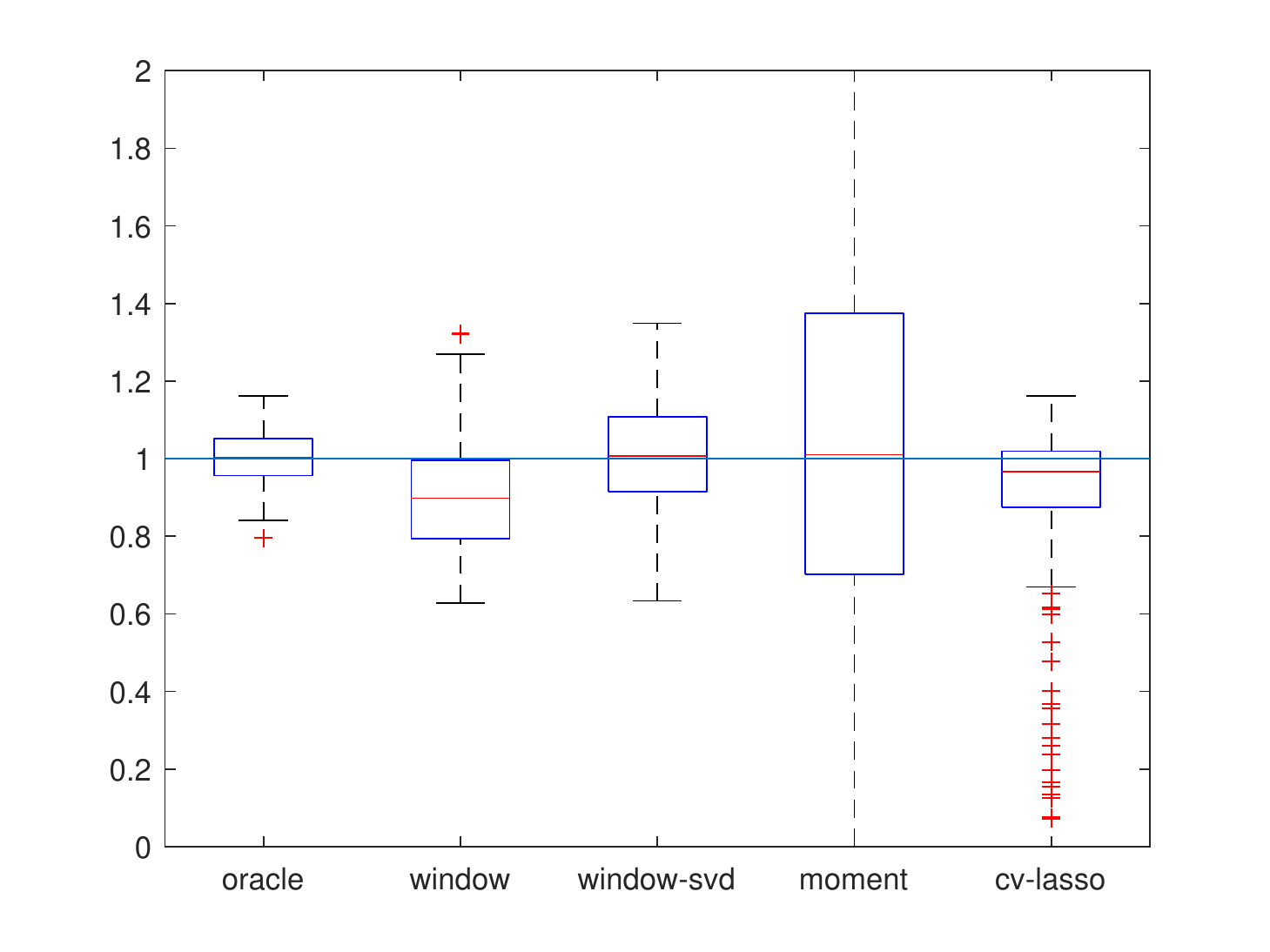}
\includegraphics[scale=.55]{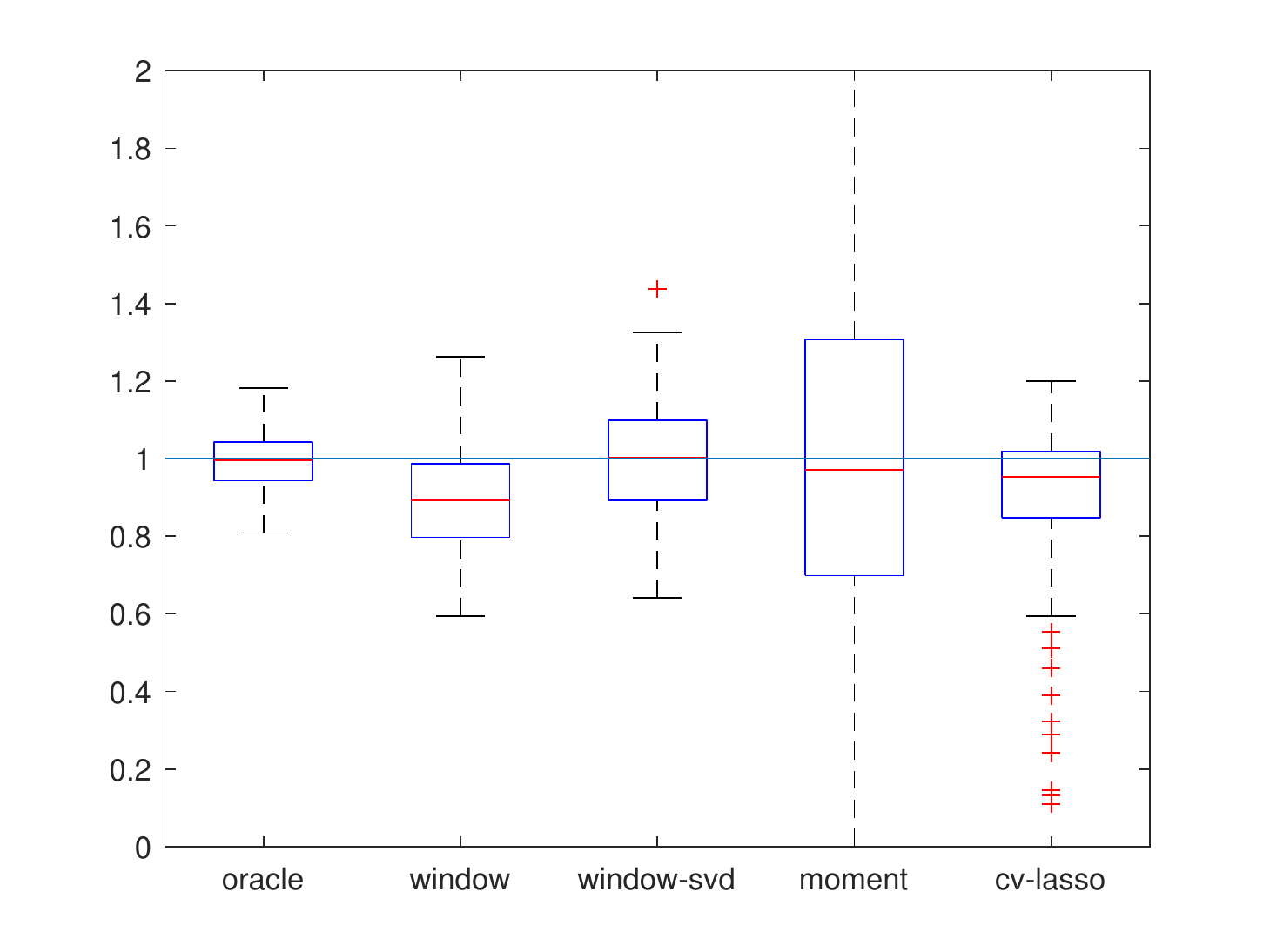}\\
  
\includegraphics[scale=.55]{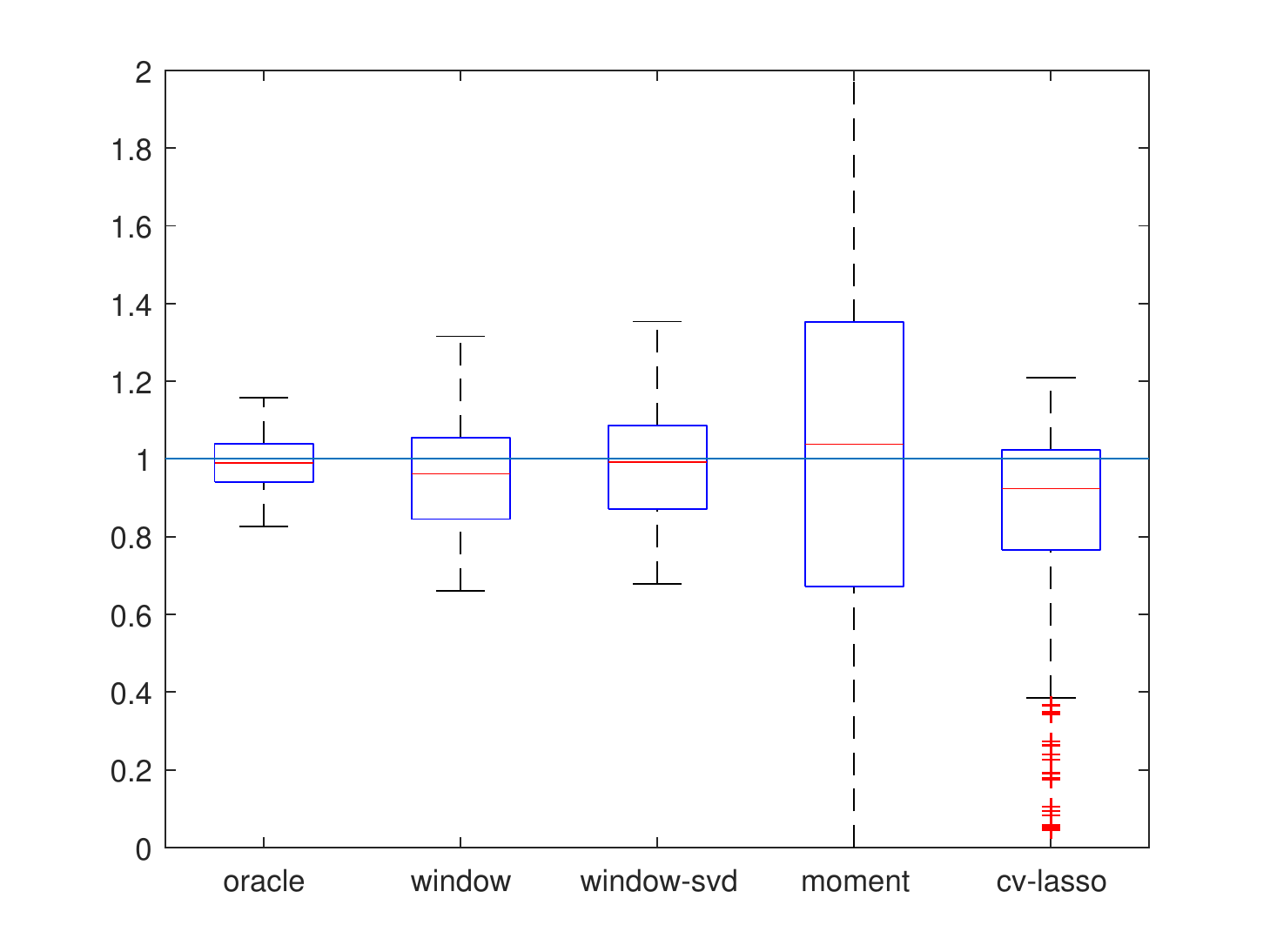}
\includegraphics[scale=.55]{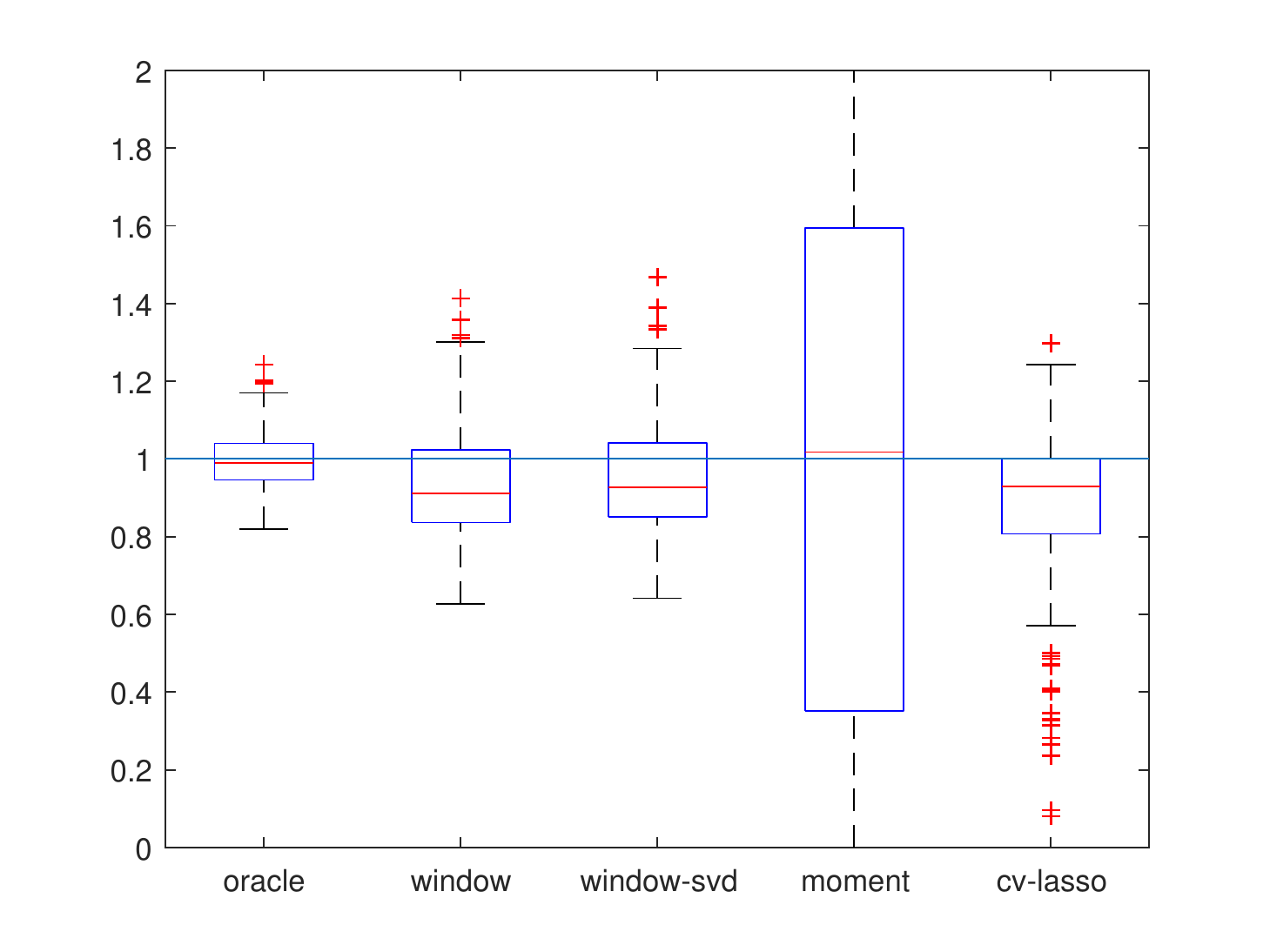}
    
    \caption{Optimal window size, $p=1,000$.  Left to right, top to bottom: Signal-less ($\|\beta\|=0$), low SNR ($\alpha = 0.1$, $\|\beta\| = 1$), medium SNR ($\alpha = 0.1$, $\|\beta\| = 5$), high SNR ($\alpha = 0.1$, $\|\beta\| = 10$) respectively.}
    \label{fig:opt}
\end{figure}
\begin{figure}
    \centering

  \includegraphics[scale=.55]{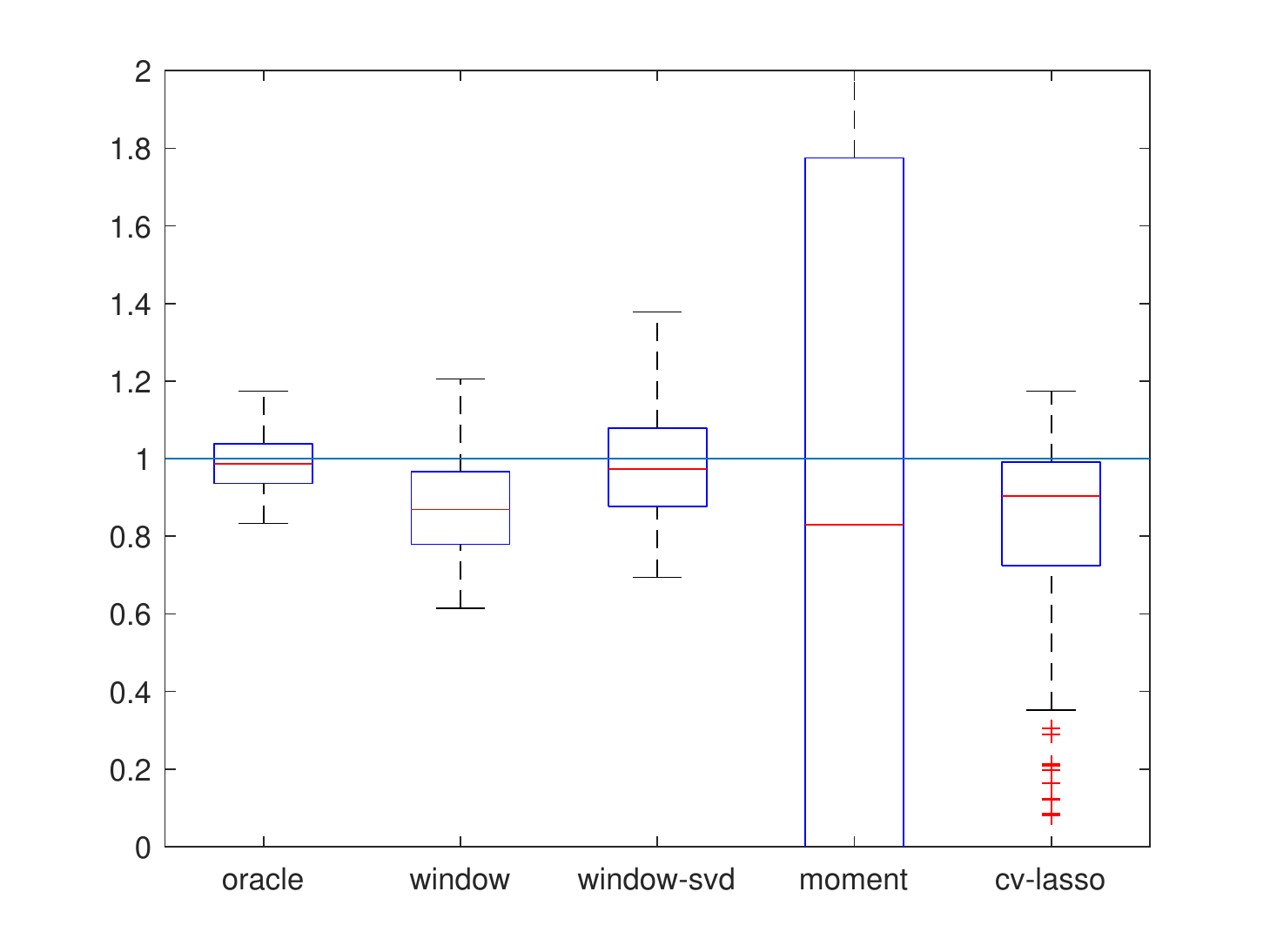}
  \includegraphics[scale=.55]{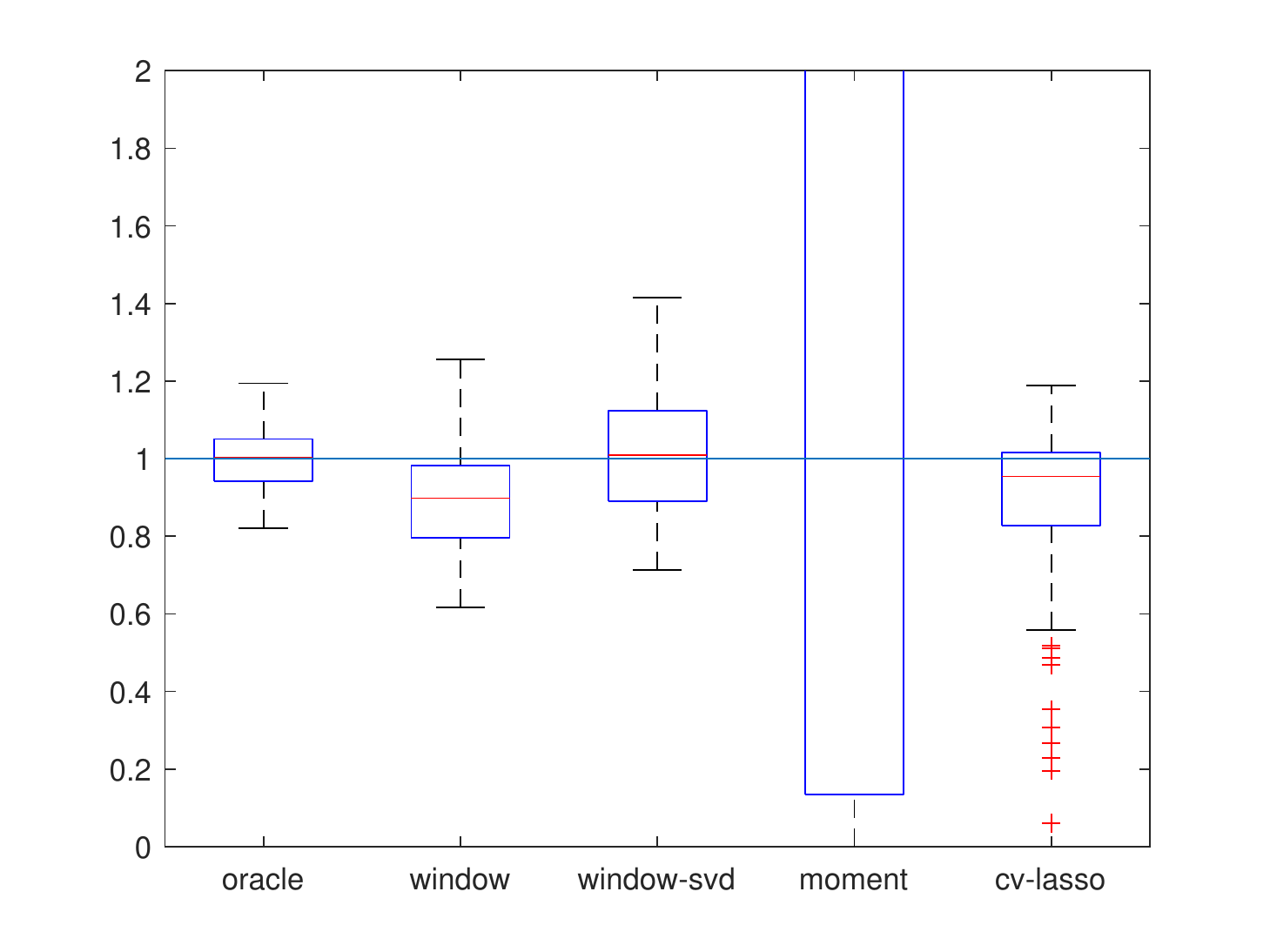}\\
  
    \includegraphics[scale=.55]{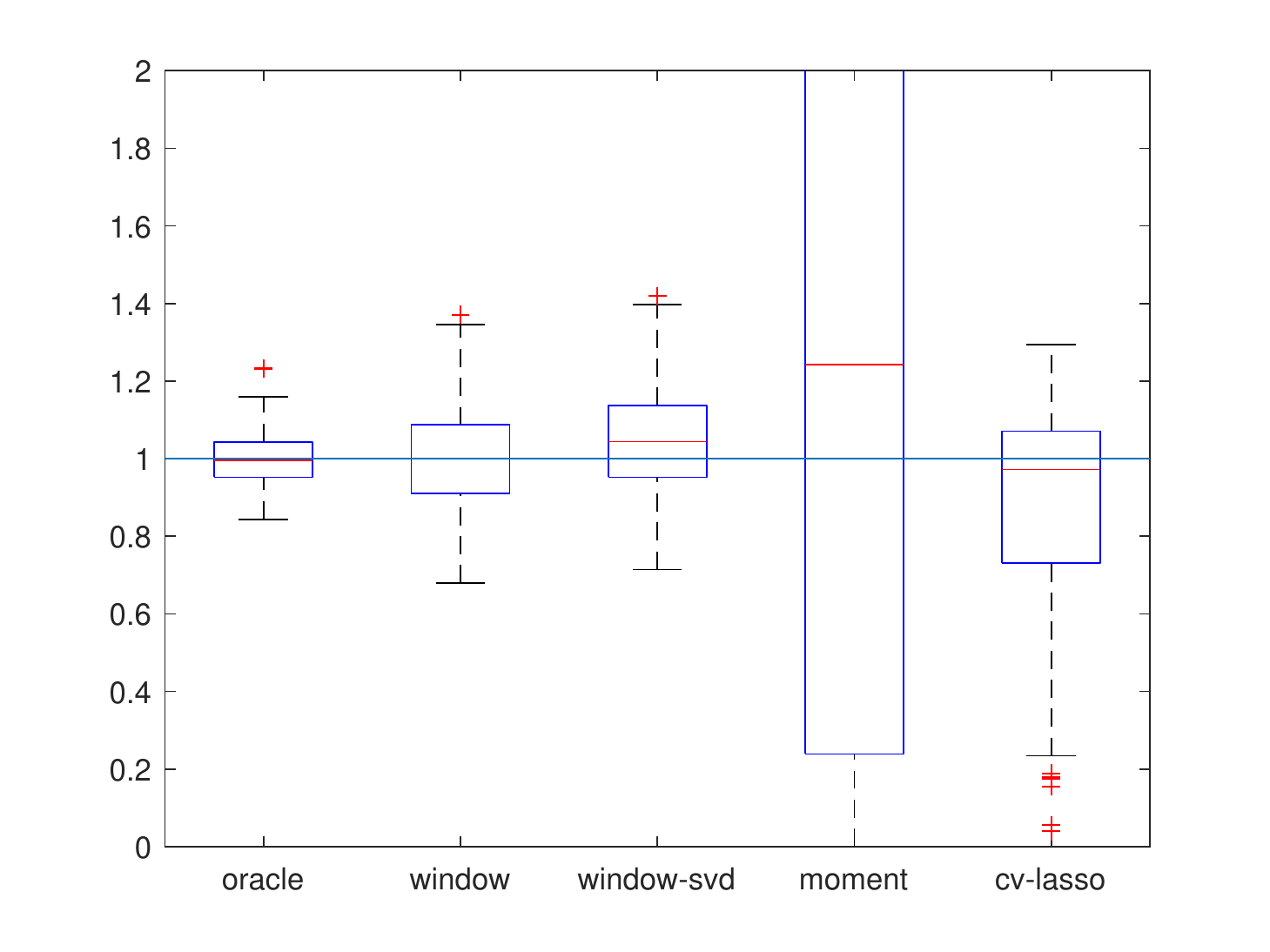}
     \includegraphics[scale=.55]{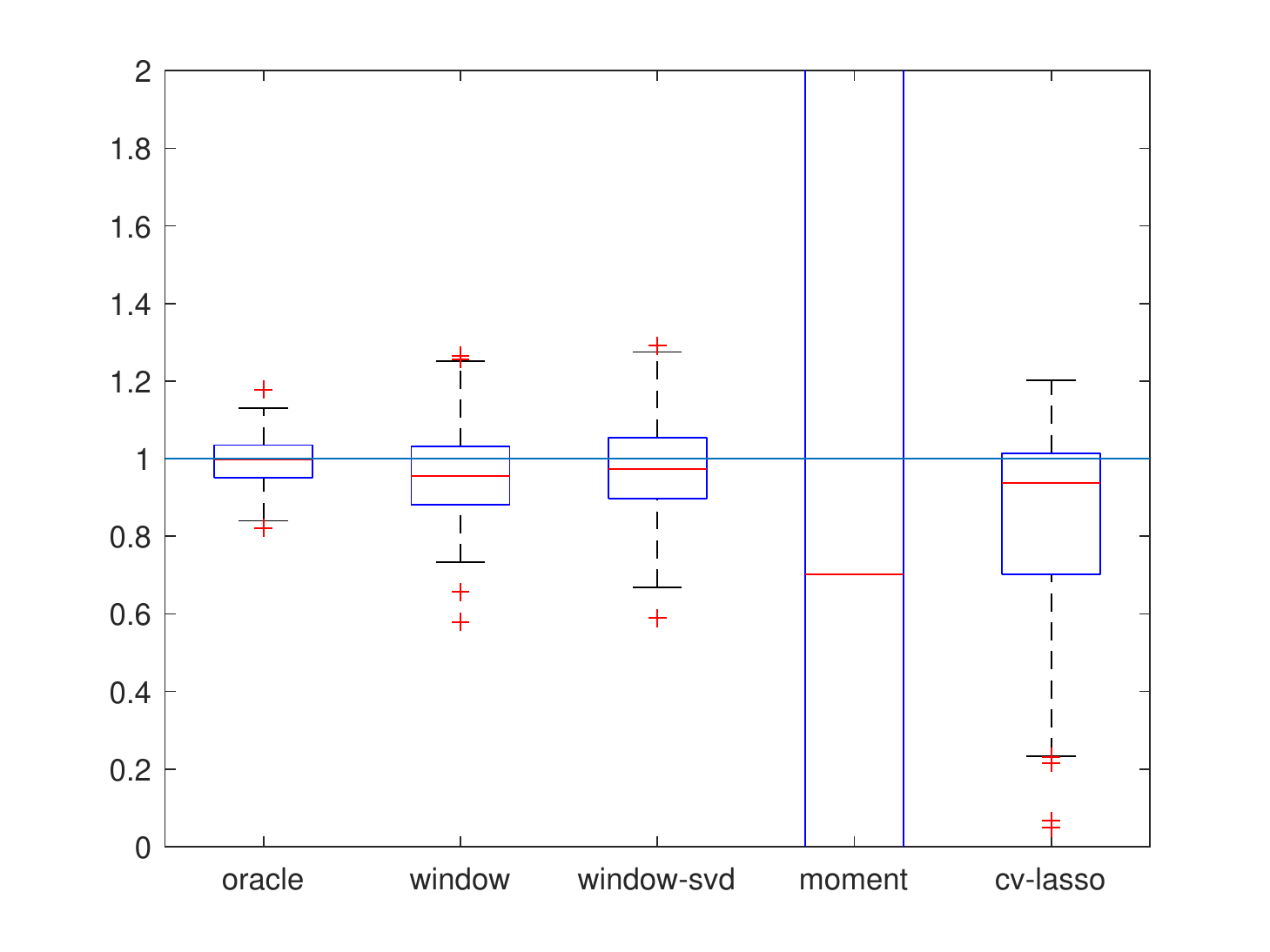}
    
    \caption{Optimal window size, $p=10,000$.  Left to right, top to bottom: Signal-less ($\|\beta\|=0$), low SNR ($\alpha = 0.1$, $\|\beta\| = 1$), medium SNR ($\alpha = 0.1$, $\|\beta\| = 5$), high SNR ($\alpha = 0.1$, $\|\beta\| = 10$) respectively.}
    \label{fig:opt_10000}
\end{figure}

\begin{figure}
    \centering

    \includegraphics[scale=.48]{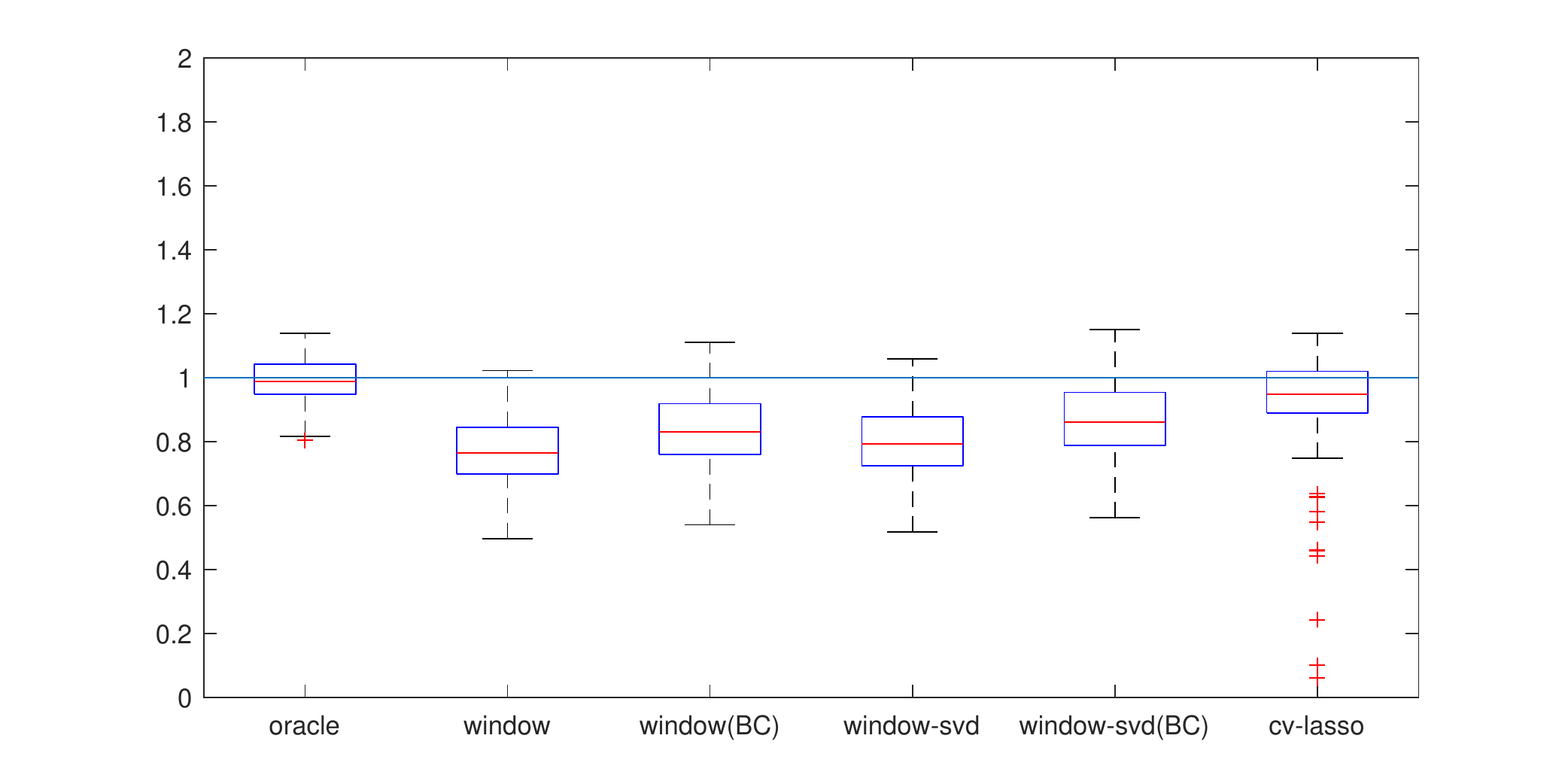}
     \includegraphics[scale=.48]{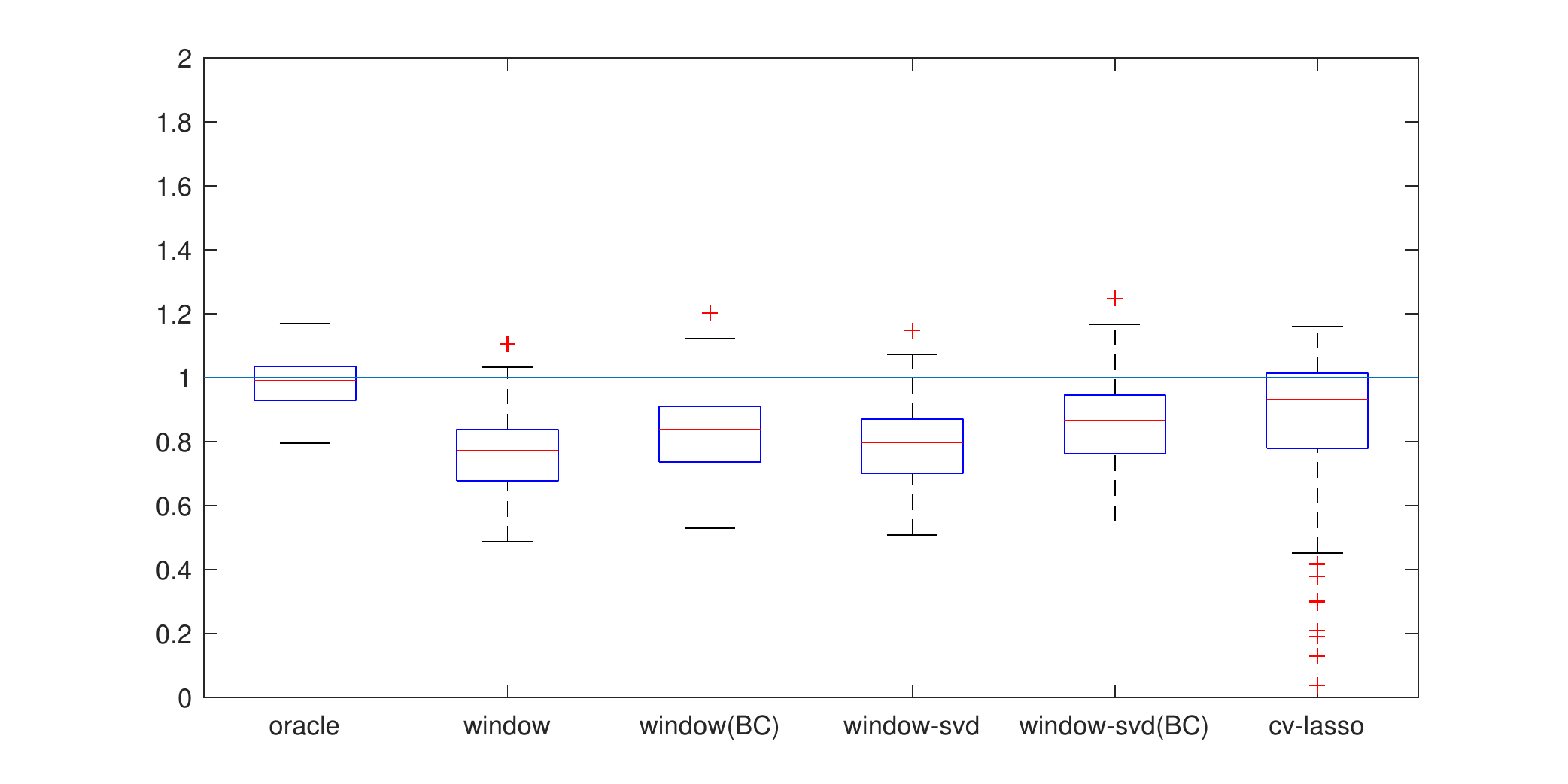}\\

    \includegraphics[scale=.48]{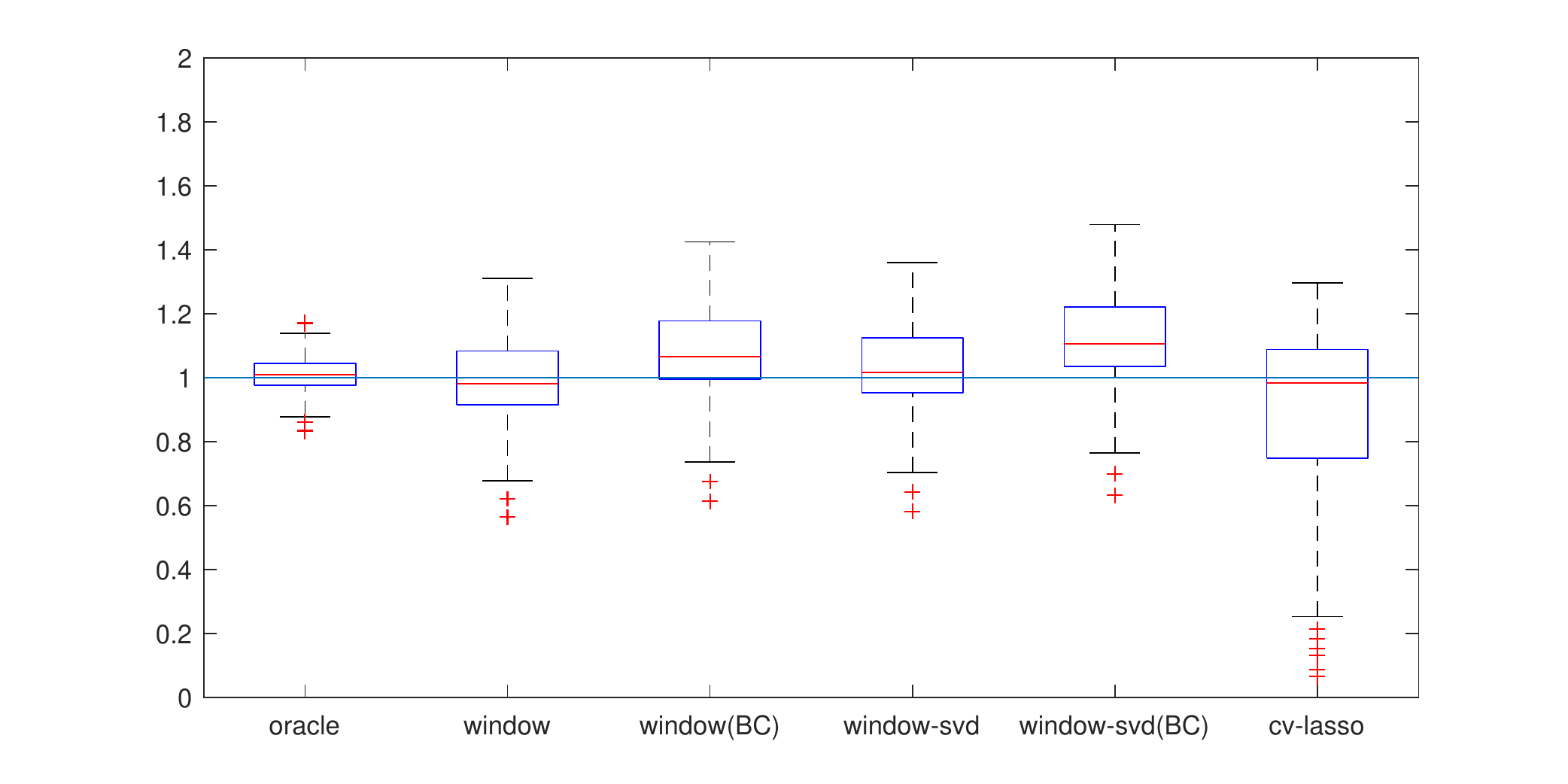}
     \includegraphics[scale=.48]{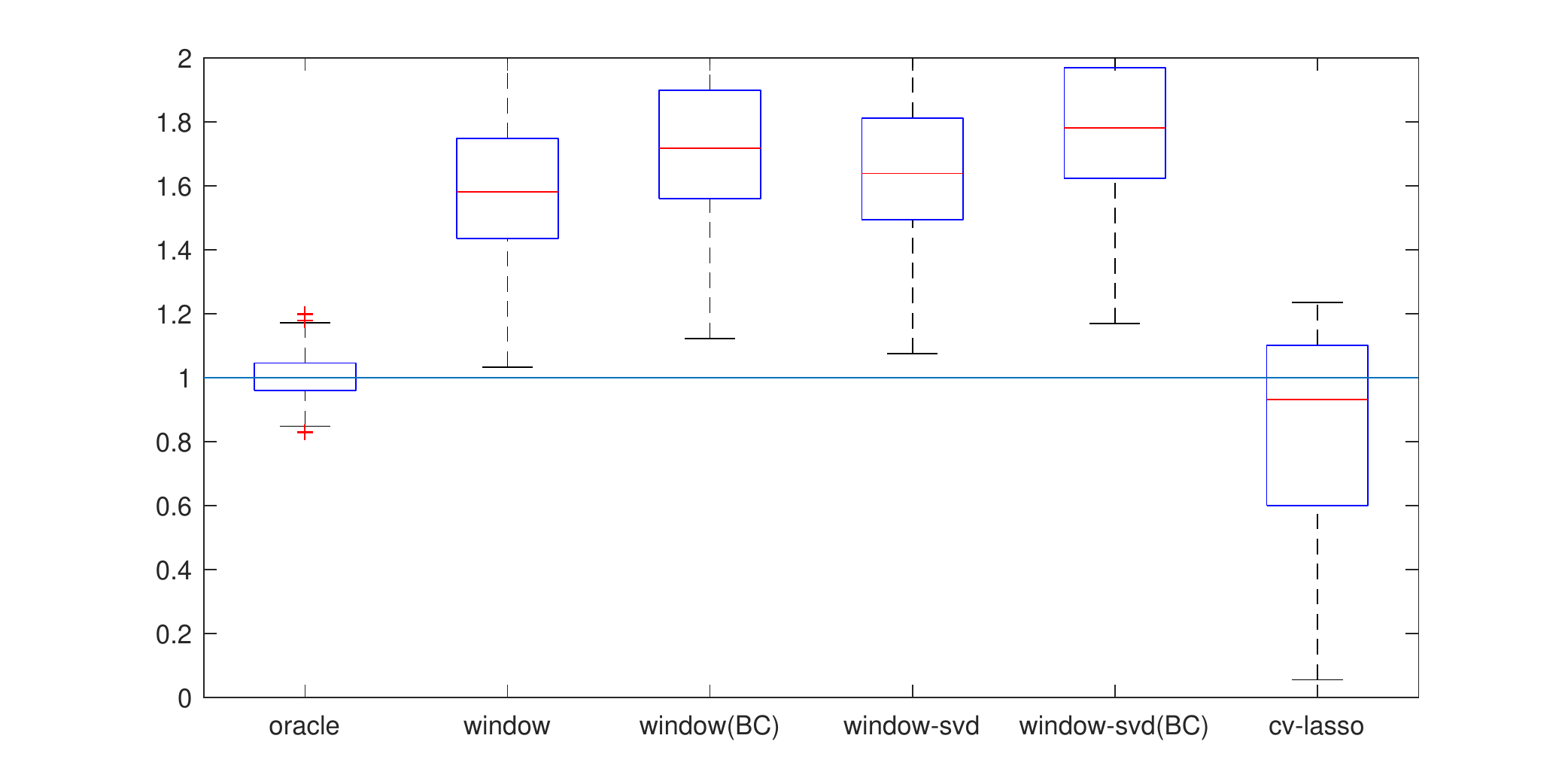}

    \caption{High-dimensional regime, $p=100000$ and (sub-optimal window size) $L=25$.   Top to bottom:  Signal-less ($\|\beta\|=0$), low SNR ($\alpha = 0.1$, $\|\beta\| = 1$), medium SNR ($\alpha = 0.1$, $\|\beta\| = 5$), high SNR ($\alpha = 0.1$, $\|\beta\| = 10$) respectively.}
    \label{fig:high}
\end{figure}

\begin{figure}
    \centering

  \includegraphics[scale=.55]{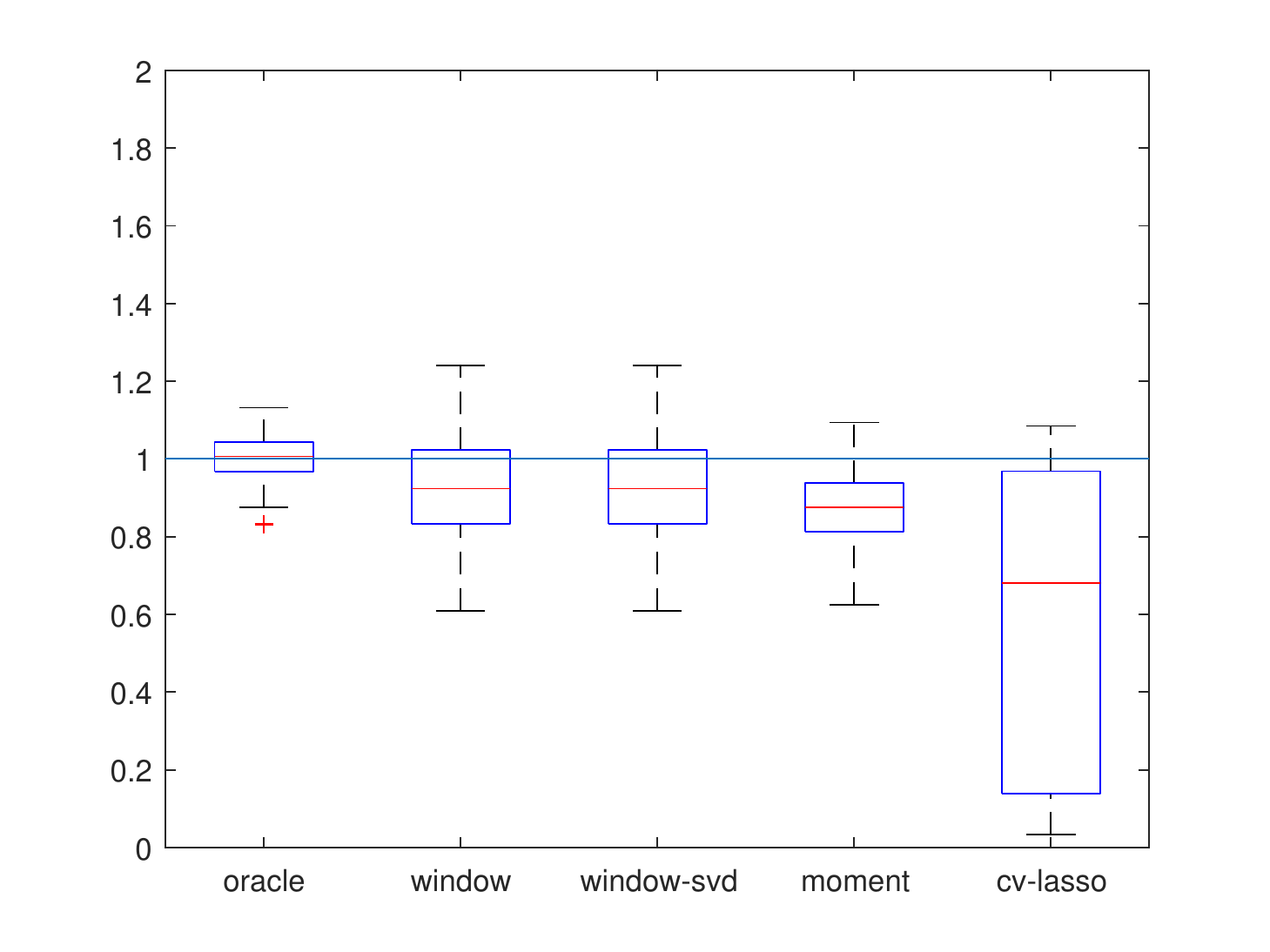} 
  \includegraphics[scale=.55]{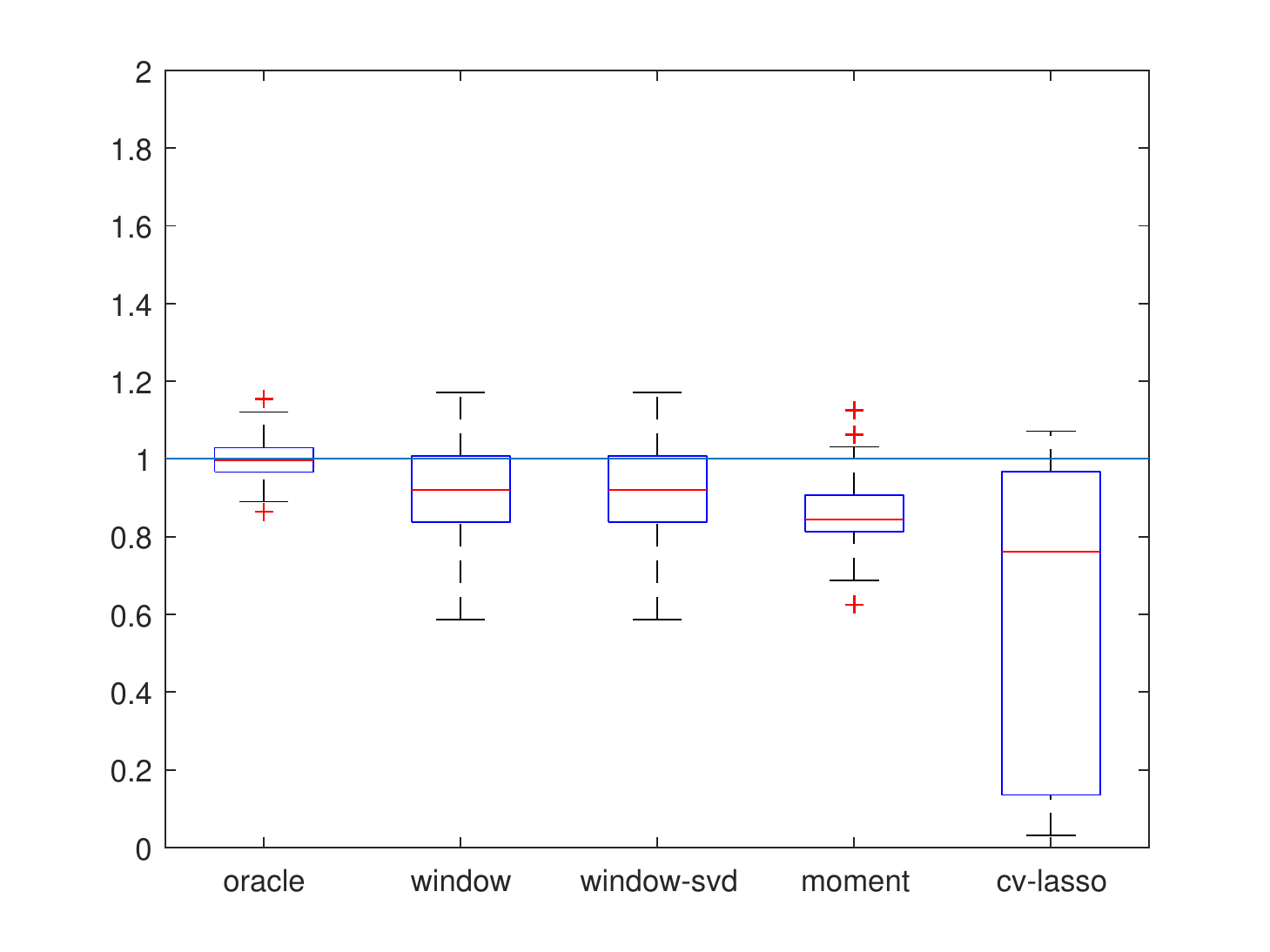}\\
  
    \includegraphics[scale=.55]{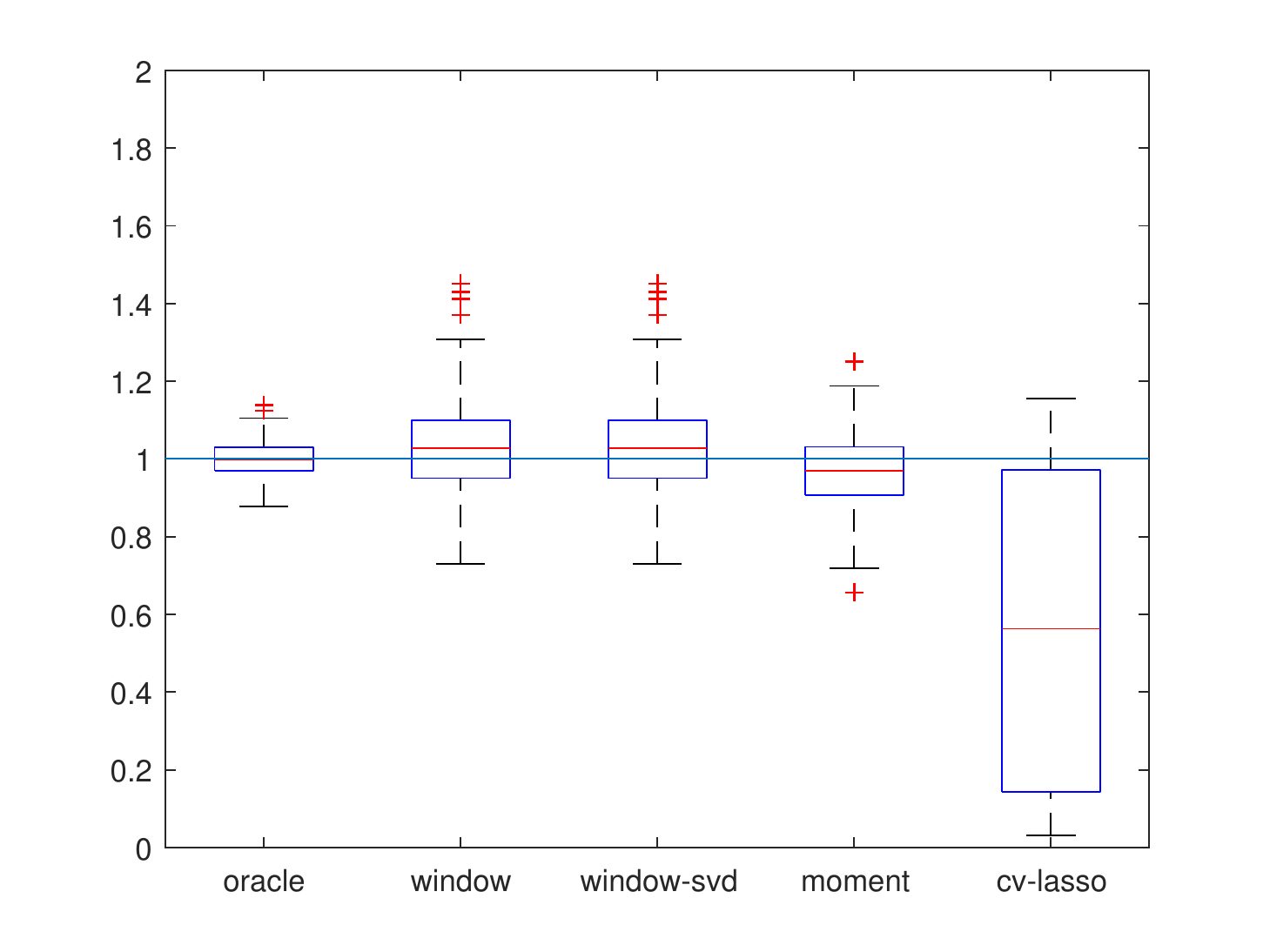}
     \includegraphics[scale=.55]{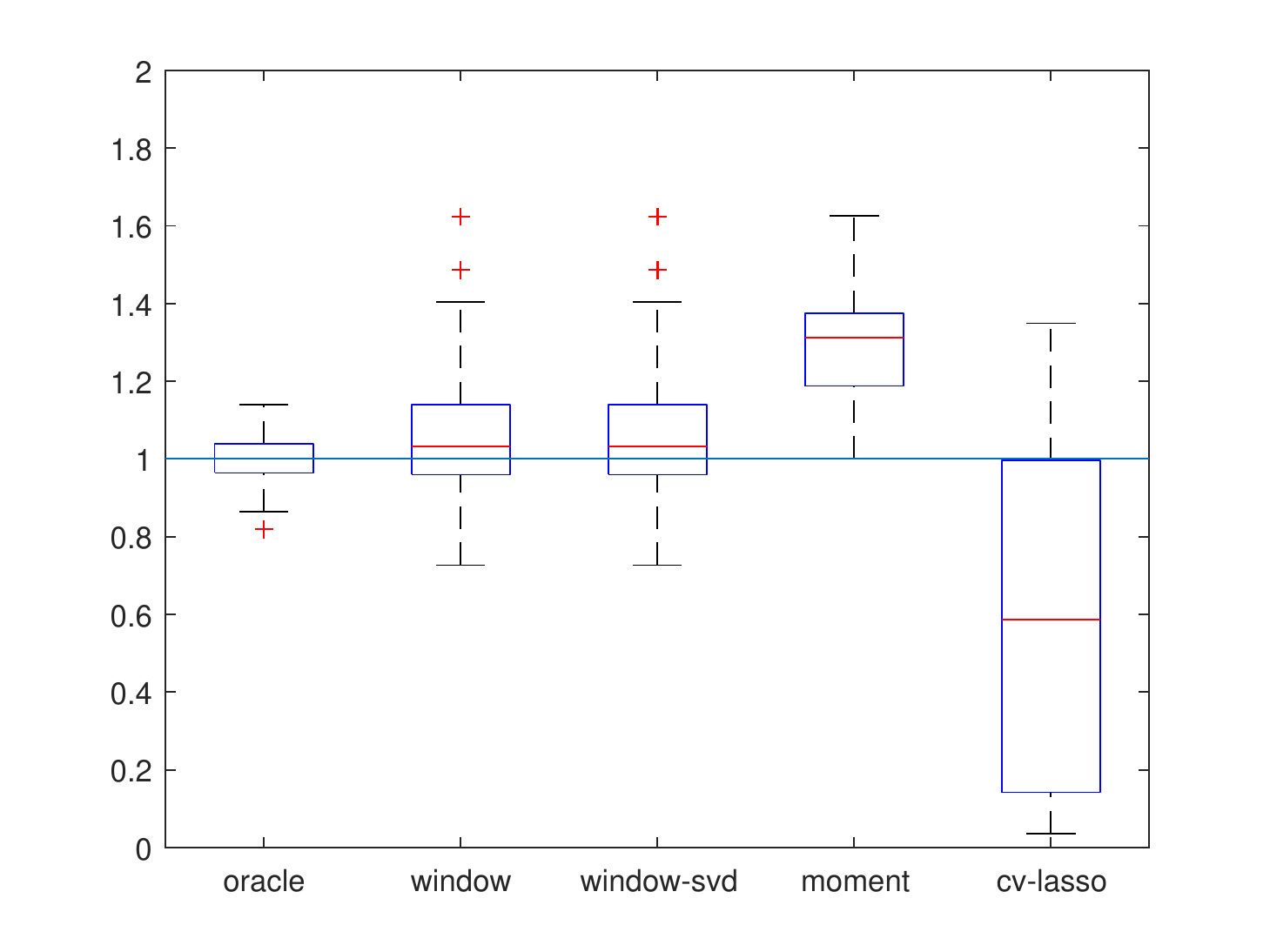}
    
    \caption{Orthogonal design matrix, $p=200$.  Left to right, top to bottom:  Signal-less ($\|\beta\|=0$), low SNR ($\alpha = 0.1$, $\|\beta\| = 1$), medium SNR ($\alpha = 0.1$, $\|\beta\| = 5$), high SNR ($\alpha = 0.1$, $\|\beta\| = 10$) respectively.}
    
    \label{fig:ortho}
\end{figure}

\begin{figure}
    \centering

  \includegraphics[scale=.55]{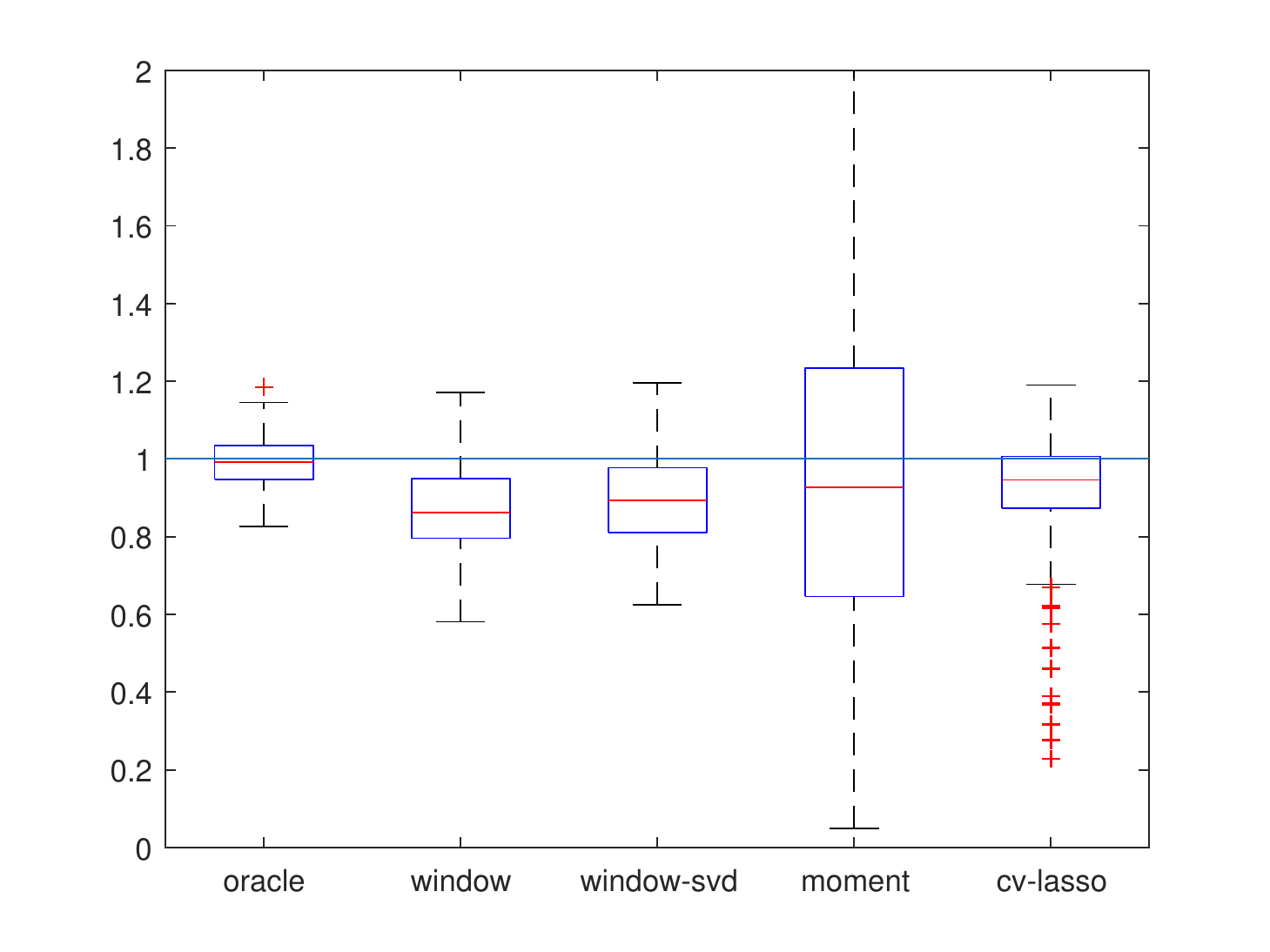}\\
  
    \includegraphics[scale=.55]{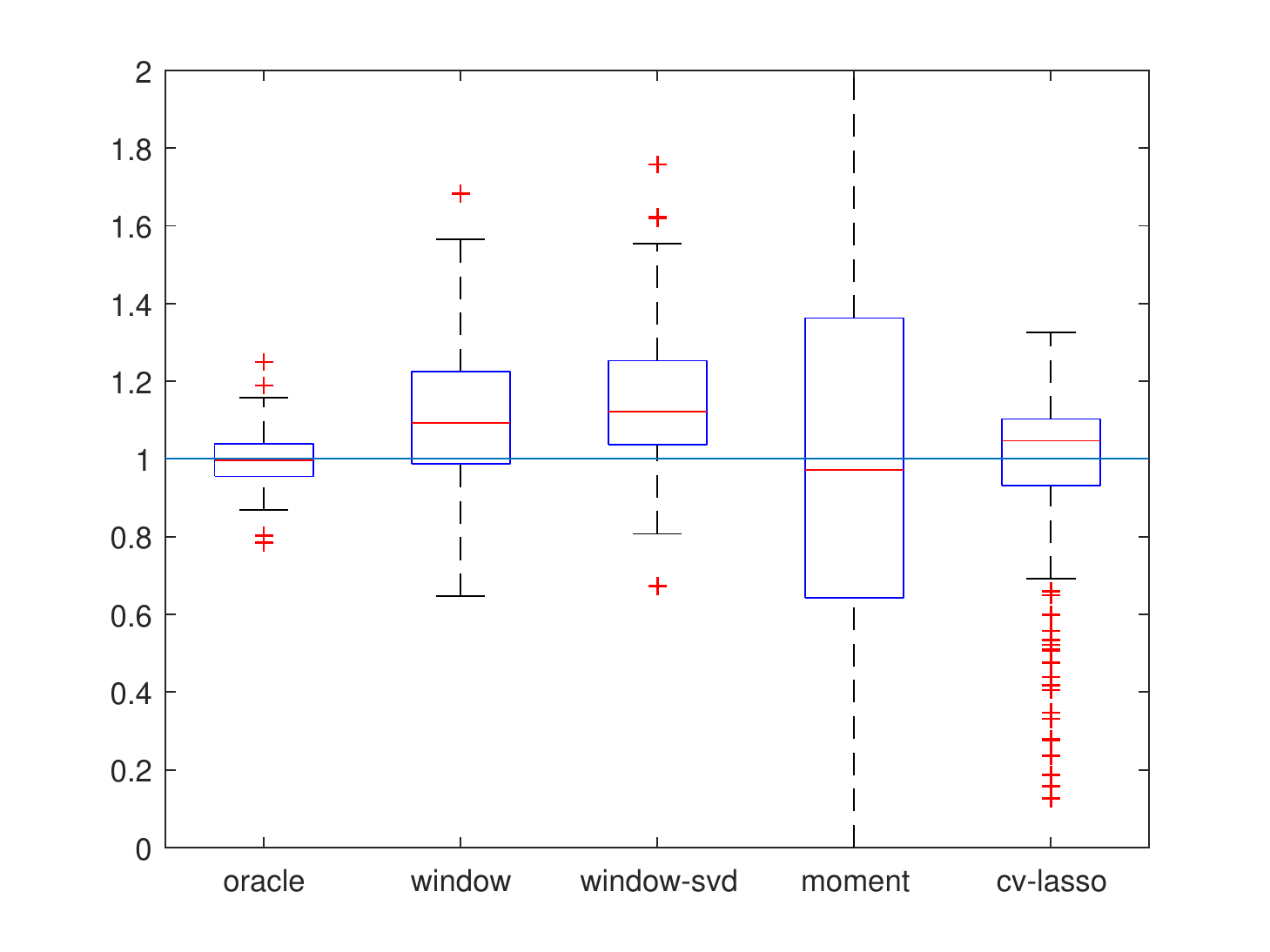}
     \includegraphics[scale=.55]{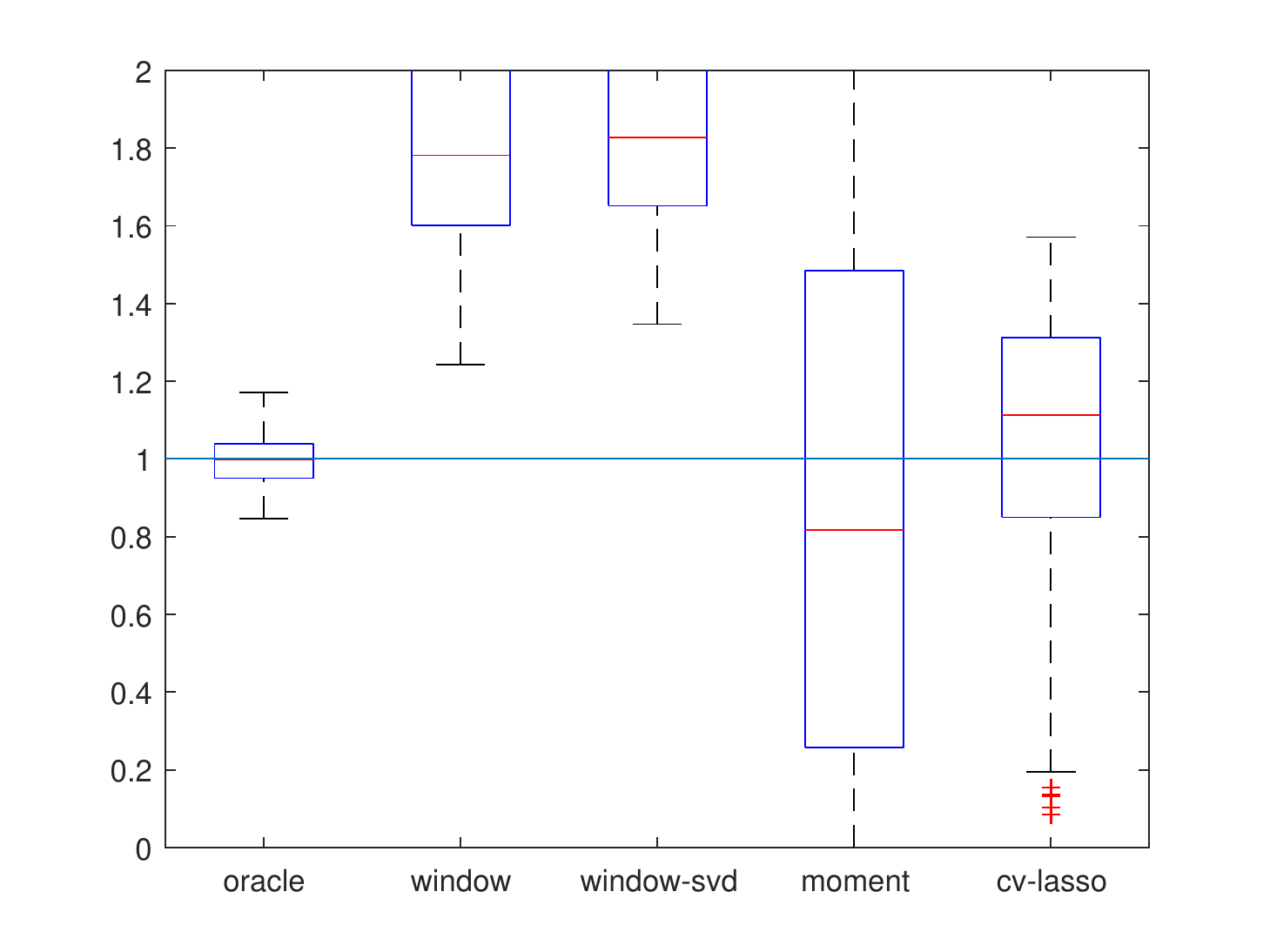}
    
    \caption{Dense signal, $p=1000$.  Left to right, top to bottom:  Low SNR ($\alpha = 0.9$, $\|\beta\| = 1$), medium SNR ($\alpha = 0.9$, $\|\beta\| = 5$), high SNR ($\alpha = 0.9$, $\|\beta\| = 10$), respectively.}
    
    \label{fig:dense}
\end{figure}

\newpage
\begin{figure}
    \centering

  \includegraphics[scale=0.85]{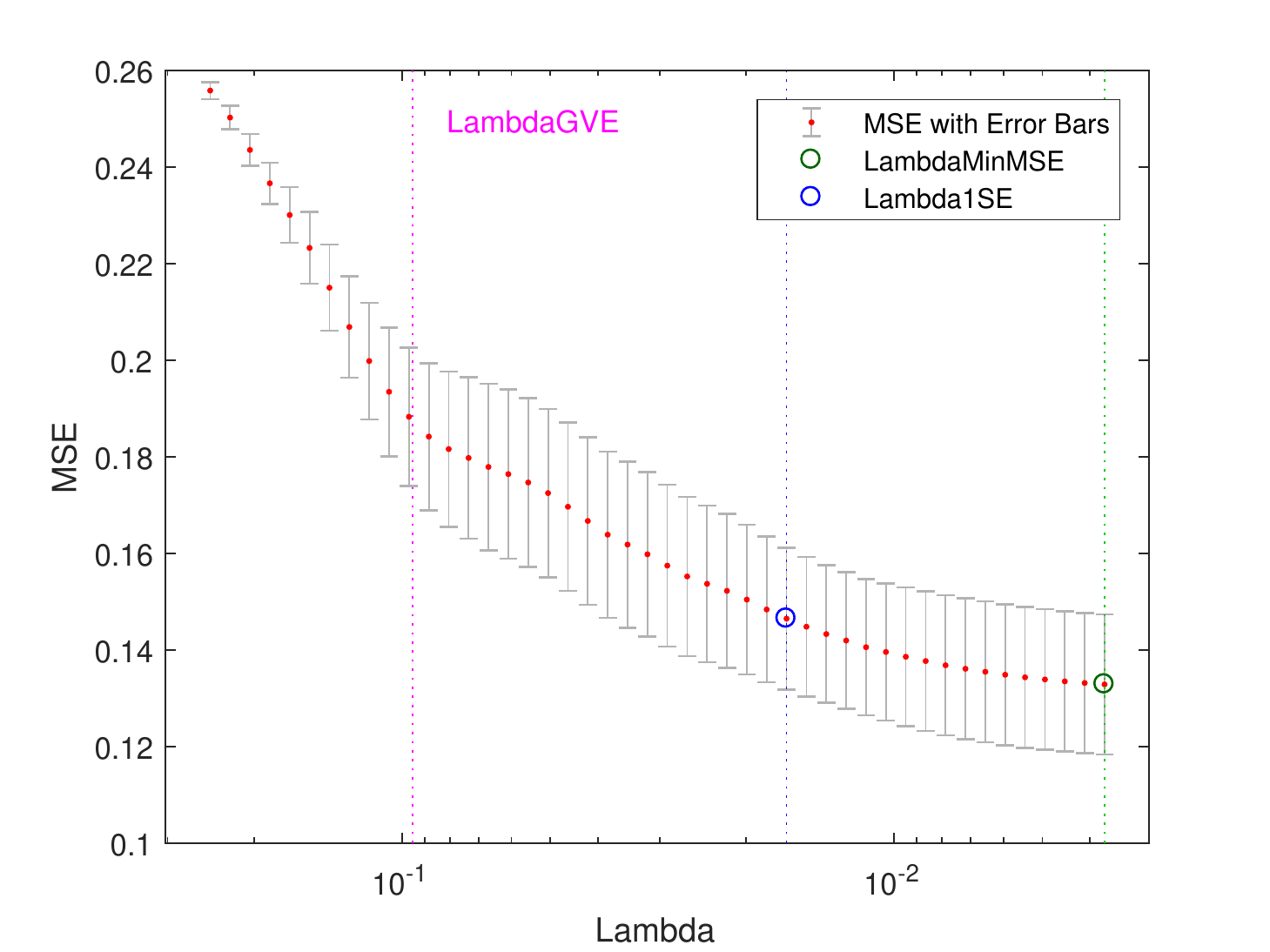}
    
    \caption{MSE for 10-fold CV LASSO using data from \cite{singh2002gene}, with the greedy variance estimator $\lambda$  from Algorithm 2 marked in magenta.}
    
    \label{fig:singh}
\end{figure}

\begin{figure}
    \centering

  \includegraphics[scale=0.85]{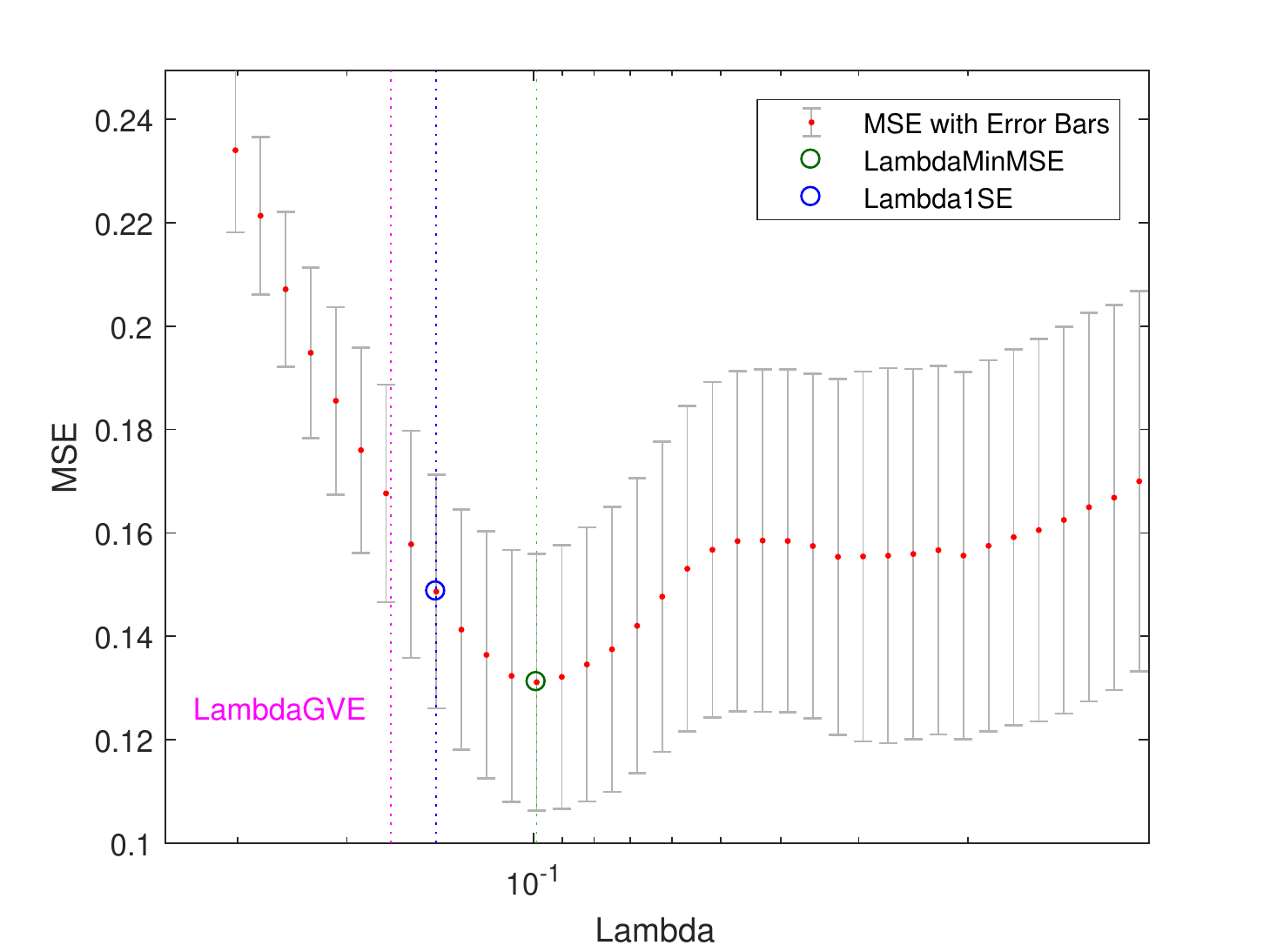}
    
    \caption{MSE for 10-fold CV LASSO using data from \cite{alon1999broad}, with the greedy variance estimator $\lambda$  from Algorithm 2 marked in magenta.}

    \label{fig:colon}
\end{figure}

\begin{figure}
    \centering

  \includegraphics[scale=0.85]{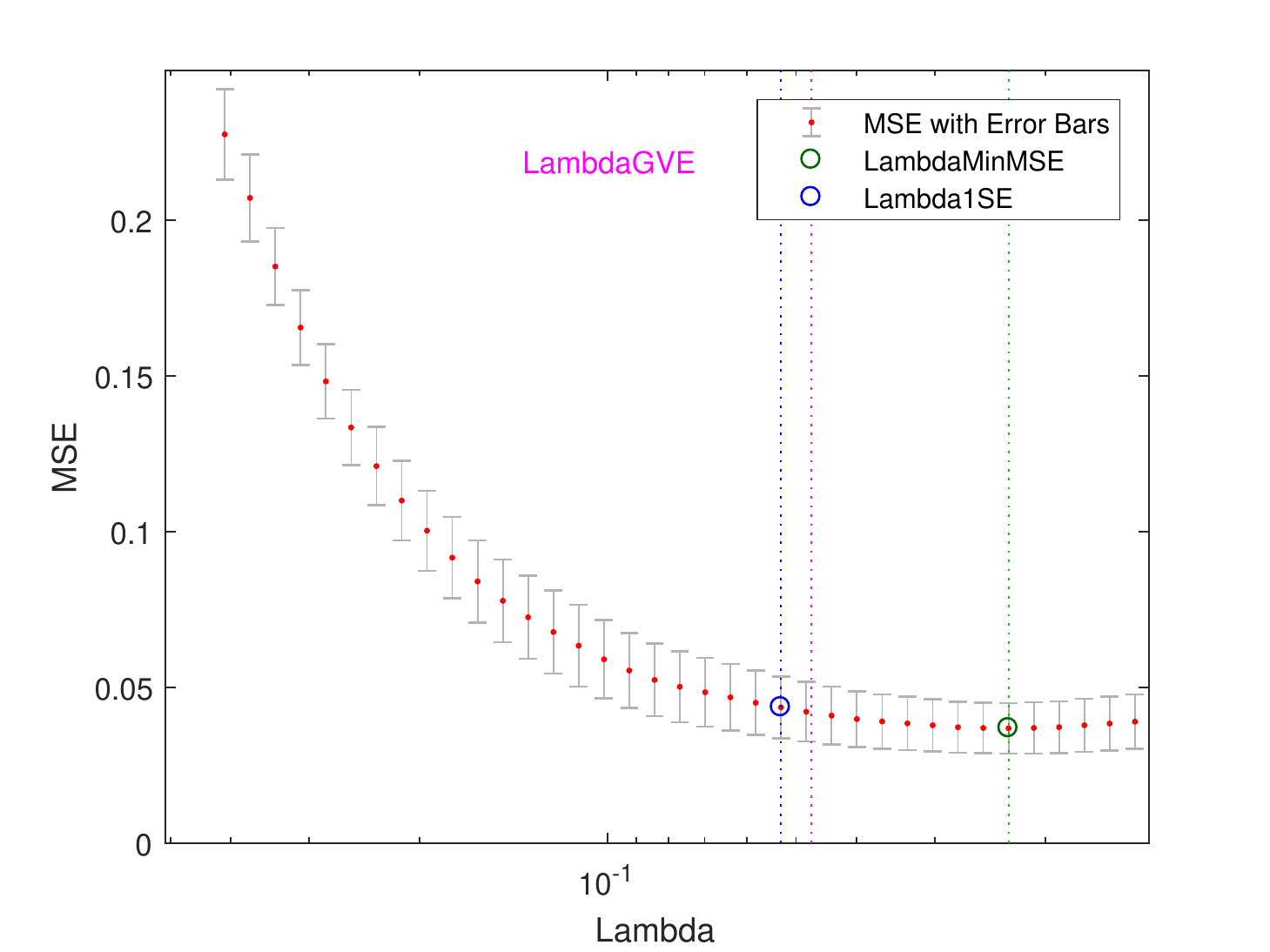}
    
    \caption{MSE for 10-fold CV LASSO using data from \cite{golub1999molecular}, with the greedy variance estimator $\lambda$  from Algorithm 2 marked in magenta.}

    \label{fig:leukemia}
\end{figure}

\end{document}